\documentclass[sort&compress,3p]{elsarticle}

\usepackage{epsfig,graphicx,subfigure}
\usepackage{makeidx}
\usepackage{multicol,multirow}
\usepackage{amsmath,amsfonts,amssymb,mathtools}
\usepackage{ulem}
\usepackage{crop,bbm}
\usepackage{xcolor,color}
\usepackage{todonotes}
\usepackage{tikz}
\usepackage{tikz-qtree}
\usepackage{subfiles}
\usepackage{stmaryrd}
\usepackage{pdflscape}
\usepackage{breqn}
\usepackage{rotating}
\usepackage{hyperref,cleveref}

\usepackage{scalerel}
\usepackage{cancel}
\usepackage{enumitem}

\usepackage{array}
\newcolumntype{C}[1]{>{\centering\arraybackslash}m{#1}}  
\newcolumntype{L}[1]{>{\raggedleft\arraybackslash}m{#1}} 
\newcolumntype{R}[1]{>{\raggedleft\arraybackslash}m{#1}} 

\newcommand{\nvar}{{m}}
\newcommand{\Id}{{\mathbf{I}}}       
\newcommand{\Idvar}{{\Id_{\nvar}}} 
\newcommand{\Idstage}{{\Id_{s}}}   
\newcommand{\Zero}{{\mathbf{0}}}  

\newcommand{\impl}{{\scaleto{\mathbf{I}}{4pt}}}
\newcommand{\expl}{{\scaleto{\mathbf{E}}{4pt}}}
\newcommand{\I}{^{\{ \impl \}}}
\newcommand{\E}{^{\{ \expl \}}}

\newcommand{\fun}{\mathbf{f}}
\newcommand{\gun}{\mathbf{g}}

\newcommand{\funim}{\fun\I}
\newcommand{\fim}{\funim}

\newcommand{\funex}{\fun\E}

\newcommand{\fex}{\funex}

 \newcommand{\xx}{\mathbf{x}}
 \newcommand{\yy}{\mathbf{y}}
 \newcommand{\zz}{\mathbf{z}}

\newcommand{\yglm}{{\mathbbm{y}}}
\newcommand{\xglm}{\mathbbm{x}}
\newcommand{\zglm}{\mathbbm{z}}

\newcommand{\Nordsieck}{\boldsymbol{\eta}}

\renewcommand{\Re}{\mathbb{R}}
\newcommand{\Co}{\mathbb{C}}
\newcommand{\one}{\mathbf{1}}

\newcommand{\Vertiii}{\vert\kern-0.22ex\vert\kern-0.22ex\vert}
\newcommand{\Vertiv}{\vert\kern-0.22ex\vert\kern-0.22ex\vert\kern-0.22ex\vert}

\newtheorem{theorem}{Theorem}
\newtheorem{lemma}{Lemma}
\newtheorem{corollary}{Corollary}
\newtheorem{remark}{Remark}
\newtheorem{definition}{Definition}
\newtheorem{assumption}{Assumption}
\newproof{proof}{Proof}
\newcommand{\norm}[1]{\left\Vert #1 \right\Vert}

\definecolor{lightblue}{rgb}{0,0,1}
\definecolor{lightred}{rgb}{1,0,0}

\usepackage{csl}

\journal{Applied Numerical Mathematics}

\newcommand{\titi}{Convergence Results for Implicit--Explicit General Linear Methods}

\begin{document}

\cslauthor{Adrian Sandu}
\csltitle{\titi}
\cslrevision{7}
\csltitlepage
	
\newpage	
\setcounter{page}{1}
\begin{frontmatter}
\title{\titi}
\author[aff1]{Adrian Sandu}
\ead{sandu@cs.vt.edu}
\address[aff1]{Computational Science Laboratory, Department of Computer Science, Virginia Tech, Blacksburg, VA 24060}

\begin{abstract}
This paper studies fixed-step convergence of implicit-explicit general linear methods. We focus on a subclass of schemes that is internally consistent, has high stage order, and favorable stability properties. Classical, index-1 differential algebraic equation, and singular perturbation convergence analyses results are given. For all these problems IMEX GLMs from the class of interest converge with the full theoretical orders under general assumptions. The convergence results require the time steps to be sufficiently small, with upper bounds that are independent on the stiffness of the problem.
\end{abstract}


\end{frontmatter}


\section{Introduction}

Consider the initial value problem for an autonomous system of differential equations in the form
\begin{equation}
\label{eqn:ivp}
 \yy'(t) =\fun(\yy)\,, \quad t_0 \le t \le t_F\,, \quad \yy(t_0) = \yy_0 \,,
\end{equation}
with $\yy(t) \in \Re^\nvar$ and $\fun: \Re^\nvar \rightarrow \Re^\nvar$  Lipschitz continuous.

General linear methods (GLMs) \cite{Butcher_1993_GLM} are numerical discretizations of \eqref{eqn:ivp} that generalize both Runge-Kutta methods (by computing multiple internal stages) and linear multistep methods (by transferring from one step to the next multiple pieces of information).
One step of the GLM applied to \eqref{eqn:ivp} reads \cite{Burrage_1980,Butcher_1996_GLMs,Butcher_2006_overview,Jackiewicz_2009_book}:
\begin{subequations}
\label{eqn:GLM} 
\begin{eqnarray}
\label{eqn:GLM-stage} 
 Y_i^{[n]} &=& h\sum_{j=1}^s a_{i,j} \, \fun\left(Y_j^{[n]}\right)+\sum_{j=1}^r u_{i,j} \, \yglm_j^{[n-1]}, \quad i=1,\dots,s,~ \\
\label{eqn:GLM-solution} 
 \yglm_i^{[n]}&=& h\sum_{j=1}^s b_{i,j} \, \fun\left(Y_j^{[n]}\right) + \sum_{j=1}^r v_{i,j} \, \yglm_j^{[n-1]}, \quad i=1,\dots,r,~
\end{eqnarray}
\end{subequations}
where $h$ is the step size. 
The method \eqref{eqn:GLM}  can be written in vector form:
\begin{subequations}
\label{eqn:GLM_vector} 
\begin{eqnarray}
\label{eqn:GLM_vector_stage} 
  Y^{[n]} &=& h \, \left(\mathbf{A} \otimes \Idvar \right) \, \fun\bigl(Y^{[n]}\bigr)+ \left(\mathbf{U} \otimes \Idvar \right)\, \yglm^{[n-1]},~ \\
\label{eqn:GLM_vector_solution} 
 \yglm^{[n]} &=& h \, \left( \mathbf{B} \otimes \Idvar \right) \, \fun\bigl(Y^{[n]}\bigr) + \left(\mathbf{V} \otimes  \Idvar \right)\, \yglm^{[n-1]},
\end{eqnarray}
\end{subequations}
where $\Idvar $ is an identity matrix of the dimension of the ODE system, the abscissa vector is $\mathbf{c} \in \Re^{s}$, and the four coefficient matrices are $\mathbf{A} \in \Re^{s \times s}$, $\mathbf{U} \in \Re^{s \times r}$, $\mathbf{B} \in \Re^{r \times s}$ and $\mathbf{V} \in \Re^{r \times r}$. The method can be represented compactly in the following Butcher tableau: 
\begin{equation}
\label{eqn:GLM-Butcher-tableau}
  \begin{array}{c|c|c}
\mathbf{c} &   \mathbf{A} & \mathbf{U} \\ \hline
  &  \mathbf{B} & \mathbf{V}
\end{array}\,.
\end{equation}
The method \eqref{eqn:GLM_vector} builds $s$ {\it internal stages} $Y_i^{[n]}$ which are meant to be order $q$ approximations of the exact solution at the abscissae:
\begin{equation}
\label{eqn:GLM-stage-order}
Y^{[n]} = \yy(t_{n-1} + \mathbf{c} h) + \mathcal{O}(h^{q+1}),
\end{equation}
and $r$ {\it external stages} $\yglm^{[n]}$ which are constructed to be order $p$ approximations of linear combinations of solution derivatives:
\begin{equation}
\label{eqn:GLM-order}
\yglm^{[n]} 
=  \left( \mathbf{W}\otimes \Idvar \right)\,\Nordsieck_p(h,\yy,t_{n}) + \mathcal{O}\left(h^{p+1}\right),
\end{equation}
where the matrix $\mathbf{W}$ and the Nordsieck vector $\Nordsieck_p$ are:
\begin{equation}
\label{eqn:Nordsieck-vector}
 \mathbf{W}=[\mathbf{w}_0 \quad \mathbf{w}_1 \quad \cdots \quad \mathbf{w}_p] \in \Re^{r \times (p+1)}, \quad
\Nordsieck_p(h,\yy,t) \coloneqq \begin{bmatrix} \yy(t)^T & \cdots & h^p\, \left( \yy^{(p)} (t)\right)^T \end{bmatrix}^T.
\end{equation}
\begin{theorem}[GLM classical order conditions \cite{Jackiewicz_2009_book}]
\label{thm:GLM-order-conditions}
A GLM \eqref{eqn:GLM} is called {\it preconsistent} if \cite{Jackiewicz_2009_book}:
\begin{equation}
\label{eqn:GLM-preconsistency}
\mathbf{U}\, \mathbf{w}_0 = \mathbf{1}_{s \times 1}\,, \quad \mathbf{V} \, \mathbf{w}_0 = \mathbf{w}_0.
\end{equation}
A preconsistent GLM \eqref{eqn:GLM} has order $p$ \eqref{eqn:GLM-order} and stage order $q=p-1$ or  $q=p$ \eqref{eqn:GLM-stage-order} if and only if:
\begin{subequations}
\label{eqn:GLM-order_conditions}
\begin{eqnarray}
\label{eqn:GLM-order_condition-1}
\frac{\mathbf{c}^{\times k}}{k!} - \frac{\mathbf{A}\,\mathbf{c}^{\times (k-1)}}{(k-1)!} - \mathbf{U} \,\mathbf{w}_k = \Zero, \quad k=1,\dots,q, \\
\label{eqn:GLM-order_condition-2}
\sum_{\ell=0}^{k}  \frac{\mathbf{w}_{k-\ell}}{\ell!} 
- \frac{\mathbf{B}\, \mathbf{c}^{\times (k-1)}}{(k-1)!} - \mathbf{V}\, \mathbf{w}_k = \Zero, \quad
 k=1,\dots,p,
\end{eqnarray}
\end{subequations}
\end{theorem}
where $\mathbf{c}^{\times k}$ denotes the component-wise power operation. Application of the GLM \eqref{eqn:GLM} to the Dahlquist model problem for error propagation $\yy' = \lambda\,\yy$
leads a numerical solution of the form
\begin{equation}
\label{eqn:stability-matrix}
\yglm^{[n]} = \mathbf{M}(z)\, \yglm^{[n-1]}, \quad \mathbf{M}(z) = \mathbf{V} + z\,\mathbf{B}\,(\Id + z\,\mathbf{A} )^{-1}\,\mathbf{U},
\end{equation}
where $z = h\lambda$ and $\mathbf{M}(z) \in \Re^{r \times r}$ is the stability matrix of the method. The stability region of the method is defined as:
\[
\mathcal{S} = \{ z \in \Co^-~:~ \sup_n \,\Vert \mathbf{M}(z)^n \Vert \le C \}.
\]

Many multiphysics systems of interest in science and engineering \eqref{eqn:ivp} are driven by stiff (fast) and non-stiff (slow) processes acting simultaneously. The implicit-explicit (IMEX) approach seeks to obtain efficient numerical solutions by solving the stiff components with an implicit, numerically stable scheme, and the non-stiff component with an explicit, computationally efficient scheme. The two schemes need to be coordinated carefully in order to obtain the desired accuracy and stability properties of the overall discretization.

IMEX schemes have been developed in the context of linear multistep methods
\cite{Hundsdorfer_2003_IMEX_DG,Crouzeix_1980,
Ruuth_1995_IMEXpattern,Seaid_2001_IMEX,Ascher_1995_IMEX_LMM},
Runge-Kutta methods 
\cite{Ascher_1997_IMEX_RK,Cooper_1980_ARK,Cooper_1983_ARK,Zhong_1996_ARK,Amitai_1998,Pareschi_2000,
Verwer_2004_IMEX_RKC,Kennedy_2003_ARK,Kennedy_2019_IMEX-ARK}, 
linearly-implicit Rosenbrock methods 
\cite{Calvo_2001_IMEX_RK,Zhong_1996_ARK}, and extrapolation methods \cite{Sandu_2010_extrapolatedIMEX}. 

The Partitioned GLM (PGLM) framework was developed by Sandu and collaborators over a series of papers \cite{Sandu_2012_ICCS-IMEX,Sandu_2012_IMEX-TSRK,Sandu_2014_IMEX-GLM,Sandu_2015_IMEX-TSRK,Zharovsky_2012_PhD}. Using this framework they constructed the first IMEX GLM schemes, and provided coefficients for methods of order two \cite{Sandu_2012_ICCS-IMEX} to six \cite{Sandu_2015_IMEX-TSRK,Sandu_2016_highOrderIMEX-GLM}.

Building on this work, Cardone at el. \cite{Sandu_2014_IMEX_GLM_Extrap} proposed an extrapolation approach to construct IMEX GLMs: start with the implicit method, and replace the non-stiff function values at each stage by extrapolated values based on the previous step. High order, stable IMEX GLMs have been constructed by Jackiewicz and co-workers \cite{Sandu_2015_Stable_IMEX-GLM,Jackiewicz_2017_IMEX-RK,Jackiewicz_2017_IMEX-construction}. Bras et al. \cite{Jackiewicz_2017_IMEX-GLM-IRKS} have constructed IMEX GLMs with inherited Runge Kutta stability property. Bras et al.  \cite{Jackiewicz_2018_IMEX-error} have studied the local truncation error for IMEX GLMs of order $q=p$ up to, and including, terms of order $\mathcal{O}(h^{p+2})$. 
Lang and Hundsdorfer used the extrapolation approach to construct IMEX GLMs of PEER type \cite{Lang_2017_IMEX-Peer}. Schneider et al. then extended the approach to superconvergent IMEX PEER schemes \cite{Schneider_2018_IMEX-Peer,Schneider_2019_IMEX-Peer}. Soleimani et al. also constructed new IMEX PEER methods of high order \cite{Weiner_2017_IMEX-Peer,Soleimani_2018_IMEX-Peer}.

This paper performs a convergence study of IMEX general linear methods for step sizes that are small enough, but do not depend on the stiffness of the problem. To the best of author's knowledge, such results were not previously available in the literature. 
\begin{enumerate}
\item We first carry out a classical convergence analysis and conclude that IMEX GLMs of high stage order converge with the theoretical order for sufficiently small step sizes; the novelty of our analysis is that the step upper bounds depend on the method coefficients and on non-stiff Lipschitz constants, but are independent of the stiffness level of the fast component. 
\item Next, we carry out a convergence analysis for IMEX GLMs applied to singular perturbation problems. This analysis follows the approach of Schneider \cite{Schneider_1993_GLM-SPP} for implicit GLMs. We show that the method has a unique numerical solution. We study the convergence for index-1 DAE problems, and then the convergence for very stiff singular perturbation problems. Under mild assumptions high stage order IMEX GLMs converge with the full theoretical order for both the non-stiff and the algebraic/stiff variables.
\end{enumerate}

The remainder of the paper is organized as follows.
A review of implicit-explicit general linear methods in the partitioned general linear method framework is provided in 
Section \ref{sec:IMEX-GLM}.
Section \ref{sec:classical-convergence} provides a classical convergence analysis for IMEX GLMs. Convergence results for IMEX GLMs applied to a singular perturbation problem are  given in Section \ref{sec:SPP-convergence}.
Section \ref{sec:discussion} discusses the findings of the paper.
\section{Implicit-explicit general linear methods for component and additively partitioned systems of differential equations}
\label{sec:IMEX-GLM}

Here we focus on partitioned GLMs \cite{Sandu_2012_ICCS-IMEX,Sandu_2012_IMEX-TSRK,Sandu_2014_IMEX-GLM,Sandu_2015_IMEX-TSRK,Sandu_2016_highOrderIMEX-GLM} applied to integrate systems with two components: a non-stiff one described by a right hand side function $\fex$, and a stiff one described by a right hand side function $\fim$. We consider partitioned GLMs where an explicit method used to solve the non-stiff component is paired with an implicit method to solve the stiff component. In this section we follow the presentation in \cite{Sandu_2014_IMEX-GLM,Sandu_2016_highOrderIMEX-GLM}.

\subsection{Component partitioned systems}
\label{sec:IMEX-component}
%
Consider the two-way component partitioned system:
\begin{subequations}
\label{eqn:component-partitioning}
\begin{eqnarray}
\label{eqn:component-partitioning-nonstiff}
\xx' &=& \fex(\xx,\zz) \,,\\
\label{eqn:component-partitioning-stiff}
\zz' &=&  \fim(\xx,\zz),
\end{eqnarray}
\end{subequations}
where $\fex$ is non-stiff and $\fim$ is stiff. For simplicity of notation we consider $\xx(t),\zz(t) \in \Re^\nvar$, although they can have different dimensions. In order to solve \eqref{eqn:component-partitioning} with an IMEX approach we start with a partitioned GLM consisting of one implicit and one explicit component:
\begin{equation}
\label{eqn:compact-imex-glm}
  \begin{array}{c|c|c|c|c}
\mathbf{c} &   \mathbf{A}\I & \mathbf{U}\I &   \mathbf{A}\E & \mathbf{U}\E \\[2pt] \hline 
  & \mathbf{B}\I & \vspace*{2pt} \mathbf{V}\I & \mathbf{B}\E & \vspace*{2pt} \mathbf{V}\E 
\end{array}.
\end{equation}
The IMEX method  \eqref{eqn:compact-imex-glm} is internally consistent in the sense that the two components share the same abscissa vector: $\mathbf{c}\I=\mathbf{c}\E=\mathbf{c}$.
Equation \eqref{eqn:component-partitioning-nonstiff} is discretized with the explicit GLM, and equation \eqref{eqn:component-partitioning-stiff} with the implicit GLM in \eqref{eqn:compact-imex-glm}, to obtain the following IMEX GLM scheme.
\begin{definition}[IMEX-GLM methods for component partitioned systems]
One step of an implicit-explicit general linear method applied to \eqref{eqn:component-partitioning} advances the solution using
\begin{subequations}
\label{eqn:component-imex-glm}
\begin{eqnarray}
\label{eqn:component-imex-glm-X}
X^{[n]} &=& h\, \left(\mathbf{A}\E \otimes \Idvar \right) \, \fex(X^{[n]},Z^{[n]}) + \left(\mathbf{U}\E \otimes \Idvar \right) \, \xglm^{[n-1]},\qquad \\ 
\label{eqn:component-imex-glm-Z}
Z^{[n]} &=& h\,\left( \mathbf{A}\I \otimes \Idvar \right) \, \fim(X^{[n]},Z^{[n]})  + \left(\mathbf{U}\I \otimes \Idvar \right) \, \zglm^{[n-1]},\qquad \\ 
\label{eqn:component-imex-glm-x}
\xglm^{[n]} &=& h\, \left(\mathbf{B}\E \otimes \Idvar \right) \, \fex(X^{[n]},Z^{[n]}) + \left(\mathbf{V}\E \otimes \Idvar \right) \, \xglm^{[n-1]}, \\
\label{eqn:component-imex-glm-z}
\zglm^{[n]} &=& h\,\left( \mathbf{B}\I \otimes \Idvar \right)  \, \fim(X^{[n]},Z^{[n]})  + \left(\mathbf{V}\I \otimes \Idvar \right) \, \zglm^{[n-1]}.\qquad
\end{eqnarray}
\end{subequations}
\end{definition}
The external stages represent {\it different} linear combinations of derivatives of the $\xx$ and $\zz$ variables \eqref{eqn:GLM-order}:
\begin{equation}
\label{eqn:component-external-stages}
\xglm^{[n]}  =  \left( \mathbf{W}\E\otimes \Idvar \right)\,\Nordsieck_p(h,\xx,t_{n}) + \mathcal{O}\left(h^{p+1}\right), \quad
\zglm^{[n]}  =  \left( \mathbf{W}\I\otimes \Idvar \right)\,\Nordsieck_p(h,\zz,t_{n}) + \mathcal{O}\left(h^{p+1}\right).
\end{equation}
%

\subsection{Additively partitioned systems}
\label{sec:IMEX-construction}
We now focus on systems with a two-way additively partitioned right hand side:
\begin{equation} 
\label{eqn:imex-ode}
 \yy' = \fex(\yy) + \fim(\yy)\,,  \quad t_0 \le t \le t_F\,, \quad   \yy(t_0)=\yy_0 \in \Re^{\nvar}\,,
\end{equation}
where $\fex$ is non-stiff and $\fim$ is stiff. 
The numerical solution is obtained as follows.

\begin{definition}[IMEX-GLM methods for additively partitioned systems]
One step of an implicit-explicit general linear method applied to \eqref{eqn:imex-ode} advances the solution using
\begin{subequations}
\label{eqn:imex_glm}
\begin{eqnarray}
\label{eqn:imex_glm-Y}
Y^{[n]} &=& h\, \left(\mathbf{A}\E \otimes \Idvar \right) \, \fex(Y^{[n]}) + h\,\left( \mathbf{A}\I \otimes \Idvar \right) \, \fim(Y^{[n]})  + \left(\mathbf{U} \otimes \Idvar \right) \, \yglm^{[n-1]},\qquad \\ 
\label{eqn:imex_glm-y}
\yglm^{[n]} &=& h\, \left(\mathbf{B}\E \otimes \Idvar \right) \, \fex(Y^{[n]}) + h\,\left( \mathbf{B}\I \otimes \Idvar \right)  \, \fim(Y^{[n]})  + \left(\mathbf{V} \otimes \Idvar \right) \, \yglm^{[n-1]}.\qquad
\end{eqnarray}
\end{subequations}
The IMEX-GLM \eqref{eqn:imex_glm} is represented compactly by the following Butcher tableau:
\begin{equation}
\label{eqn:imex_glm_tableau}
\renewcommand{\arraystretch}{1.5}
  \begin{array}{c|c|c|c}
\mathbf{c} &   \mathbf{A}\E & \mathbf{A}\I & \mathbf{U} \\ \hline 
   & \mathbf{B}\E & \vspace{2pt} \mathbf{B}\I & \mathbf{V}
\end{array}.
\end{equation}
For additively partitioned systems we consider  families of GLMs that share the same weights, $\mathbf{V}\E = \mathbf{V}\I = \mathbf{V}$ and $\mathbf{U}\E = \mathbf{U}\I = \mathbf{U}$. In this case a single set of external stages $\yglm^{[n]}$ is sufficient to carry information between steps  \cite{Sandu_2012_ICCS-IMEX,Sandu_2012_IMEX-TSRK,Sandu_2014_IMEX-GLM,Sandu_2015_IMEX-TSRK,Sandu_2016_highOrderIMEX-GLM}.
We also focus on internally consistent methods where $\mathbf{c}\I=\mathbf{c}\E=\mathbf{c}$.
\end{definition}

\begin{remark}[The external stage vectors]
The combined external stages \eqref{eqn:imex_glm-y} add together {\it different} combinations of the derivatives of the two components:
\begin{equation}
\label{eqn:mixed-derivatives-in-external-stages}
 \yglm^{[n]} = \mathbf{w}_{0}\, \yy(t_n) + \sum_{k=1}^p \mathbf{w}\E_{k}\, h^k\, (\fex)^{(k-1)} (\yy(t_n)) + \sum_{k=1}^p \mathbf{w}\I_{k}\, h^k\, (\fim)^{(k-1)} (\yy(t_n))  + \mathcal{O}\left(h^{p+1}\right).
\end{equation}
The starting procedure needs to form the combined external stages \eqref{eqn:mixed-derivatives-in-external-stages} by computing separately the derivatives of $\xx$ and $\zz$. Afterwards, the combined external stages are advanced at each step as regular GLMs do. More details on starting procedures for GLMs can be found in \cite{Jackiewicz_2017_starting,Jackiewicz_2009_book,Butcher_1996_GLMs}.
\end{remark}

\subsection{Class of methods of interest}
\label{sec:interesting-methods}
%
The following fundamental property shows that the partitioned GLM framework is very well suited to construct multimethods.
\begin{theorem}[Order conditions for partitioned GLMs  \cite{Sandu_2012_ICCS-IMEX,Sandu_2012_IMEX-TSRK,Sandu_2014_IMEX-GLM,Sandu_2015_IMEX-TSRK,Sandu_2016_highOrderIMEX-GLM}]
\label{thm:GLM-partitioned-order}
An internally consistent partitioned GLM \eqref{eqn:compact-imex-glm}, \eqref {eqn:imex_glm_tableau} has order  $p$ \eqref{eqn:GLM-order} and stage order $q \in \{p-1,p\}$ \eqref{eqn:GLM-stage-order} if and only if each component method $\{\mathbf{A}\E$, $\mathbf{B}\E$, $\mathbf{U}\E$, $\mathbf{V}\E\}$ and $\{\mathbf{A}\I$, $\mathbf{B}\I$, $\mathbf{U}\I$, $\mathbf{V}\I\}$ has order $p$  \eqref{eqn:GLM-order} and stage order $q$ \eqref{eqn:GLM-stage-order}.
\end{theorem}
Theorem \eqref{thm:GLM-partitioned-order} is applicable to PGLMs constructed using any number of individual methods, and applied to systems with any number of partitions. Each component method needs to independently meet its own order conditions \eqref{eqn:GLM-order_conditions}.  No additional ``coupling'' conditions are needed for the partitioned GLM (i.e., no order conditions contain coefficients from multiple component schemes).

In this work we focus on a subset of IMEX GLMs that have several favorable properties.

\begin{definition}[Class of IMEX GLMs of interest]
\label{def:interesting-methods}
The class of ``methods of interest'' consists of IMEX-GLMs with the following properties:
\begin{enumerate}[topsep=0pt,itemsep=0pt,parsep=0pt,partopsep=0pt]
\item the method is internally consistent, i.e., the abscissae of the explicit and implicit components coincide, $\mathbf{c}\E = \mathbf{c}\I = \mathbf{c}$;
\item the stage orders of the explicit and implicit components are $q\E,q\I \in \{p-1,p\}$; 
\item the implicit component has a coefficient matrix $\mathbf{A}\I$ with (strictly) positive eigenvalues; and 
\item the spectral radius of the implicit stability matrix at infinity is $\rho\big(\mathbf{M}\I(\infty) \big)<1$.
\end{enumerate}
\end{definition}

\section{Convergence analysis for additive IMEX GLMs}
\label{sec:classical-convergence}

We next study convergence for a fixed-step $h$ solution. Consider a partitioned system \eqref{eqn:imex-ode} where the nonstiff component is Lipschitz-continuous with a moderate Lipschitz constant in a vicinity of the exact solution:
\[
\Vert \fex(\yy) - \fex(\zz) \Vert \le \mathrm{L}\E\, \Vert \yy - \zz \Vert.
\]
The Lipschitz constant of the stiff component can be arbitrarily large. To avoid using it we make the following assumption.

\begin{assumption}[Separability of stiffness]
\label{ass:stiffness}
Assume that, for any $\tau \in [t_0,t_F]$, there is an interval $[\tau - \varepsilon, \tau+\varepsilon]$ such that the implicit component can be locally decomposed into a linear part and a nonlinear remainder:
\begin{equation}
\label{eqn:linear-nonlinear-split}
\fim(\yy(t)) = \mathbf{J}_\tau\,\yy(t) + \mathbf{r}\I_\tau(\yy(t)) \qquad \forall t \in (\tau - \varepsilon, \tau+\varepsilon),
\end{equation}
where the matrix $\mathbf{J}_\tau$ enjoys the following properties:
\begin{enumerate}
\item It is diagonalizable $\mathbf{J}_\tau = \mathbf{V}\, \operatorname{diag}(\lambda_i)\, \mathbf{V}^{-1}$, and all its eigenvalues have non-positive real parts, $\lambda_i \in \Co^-$. 
\item $\mathbf{J}_\tau$ captures all the stiffness of the system in a vicinity of the exact trajectory, i.e., the remaining nonlinear part $\mathbf{r}\I_\tau(\yy)$ is non-stiff, and is Lipschitz-continuous with a moderate Lipschitz constants in a vicinity of the exact solution:
\begin{equation*}
\Vert \mathbf{r}\I_\tau(\yy) - \mathbf{r}\I_\tau(\zz) \Vert \le \mathrm{L}\I_\tau\, \Vert \yy - \zz \Vert.
\end{equation*}
The implicit remainder Lipschitz constants $\mathrm{L}\I_\tau$ are of the same order of magnitude as the nonstiff constant $\mathrm{L}\E$.
%
\end{enumerate}
\end{assumption}
Choose a finite number of $\tau_i$ values such that the corresponding intervals $(\tau_i-\varepsilon_i,\tau_i+\varepsilon_i)$ cover the entire compact integration time interval $[t_0,t_F]$. Consequently, we select a finite number of subintervals and on each we have the corresponding decomposition \eqref{eqn:linear-nonlinear-split}.

\begin{lemma}[Eigenvalues of Kronecker product]
\label{lem:eigs-kronecker}
We review the result \cite[Theorem 4.2.12]{Horn_1991_book}.
Let $A \in \Re^{m \times m}$ with eigenvalues $a_i$ and the corresponding eigenvectors $u_i$, and $B \in \Re^{n \times n}$ with eigenvalues $b_j$ and the corresponding eigenvectors $v_j$. Then the eigenvalues of the Kronecker product $A \otimes B$ are of the form $a_i\,b_j$ with the corresponding eigenvectors $u_i \otimes v_j$ for $i=1,\dots,m$ and $j=1,\dots,n$. 
\end{lemma}

In order to study the convergence of the IMEX GLM we will use the linear-nonlinear decomposition \eqref{eqn:linear-nonlinear-split} of the stiff system component in each of these subintervals. Since we have a finite number of decompositions, without loss of generality, we will carry out the convergence analysis on a single subinterval and a single decomposition \eqref{eqn:linear-nonlinear-split}; the subscripts $\tau$ will be dropped.

\begin{theorem}[Convergence of IMEX GLM schemes]
Apply an order $p$ IMEX GLM scheme from the class of interest (according to Definition \ref{def:interesting-methods}) to solve the partitioned system \eqref{eqn:imex-ode}. We make Assumption \ref{ass:stiffness} for the stiff component $\fun\I$, which admits a splitting with a linear part that captures all the stiffness. The eigenvalues $h\,\lambda$ of $h\,\mathbf{J}$ fall within the region of absolute stability of the implicit component, i.e., the implicit component is linearly stable when applied to integrate the stiff component $\fim$.

The numerical solution converges with order $p$ to the exact solution for all sufficiently small step sizes $h \le h^\ast$, where the upper bound $h^\ast$ depends on the method coefficients and on the Lipschitz constants $\mathrm{L}\E$, $\mathrm{L}\I$, but is independent of the stiffness of $\fun\I$.
%
\end{theorem}

\begin{proof}

We apply an IMEX GLM scheme from the class of interest according to Definition \ref{def:interesting-methods}. Let $a_i > 0$ be the positive eigenvalues of $\mathbf{A}\I$.
From Lemma \ref{lem:eigs-kronecker} the eigenvalues of the matrix $\mathbf{A}\I \otimes h\,\mathbf{J}$ are $a_i\, (h\lambda_j)$, where $\lambda_j$ are the eigenvalues of $\mathbf{J}$; moreover, the eigenvectors of this Kronecker product are independent of $h$. This implies that the following matrix is uniformly bounded for any step size:
\[
\mathbf{S} \coloneqq \Id_{\nvar s} - \mathbf{A}\I \otimes h\,\mathbf{J}, \quad
\Vert \mathbf{S}^{-1} \Vert \le C_{\mathbf{S}}~~\forall\, h,
\]
since $\mathbf{S}^{-1}$ has eigenvalues $(1 - a_i\, (h\lambda_j))^{-1}$, which are uniformly bounded for any $h$, and its eigenvectors are independent of $h$.

The assumption that the implicit method is linearly stable on the stiff component implies that the implicit stability matrix is power bounded:
\[
\big\Vert \mathbf{M}\I (h \mathbf{J})^n \big\Vert \le C_{\mathbf{M}}\quad \forall n.
\]
Replace the exact solution into the method \eqref{eqn:imex_glm}, then subtract the numerical solution to obtain the following recurrence of global errors:
\begin{eqnarray*}
\Delta Y^{[n]} &=& h\, \left(\mathbf{A}\E \otimes \Idvar \right) \, \Delta\fex(Y^{[n]}) + h\,\left( \mathbf{A}\I \otimes \Idvar \right) \, \left(\mathbf{J}\, \Delta Y^{[n]} + \Delta\mathbf{r}\I(Y^{[n]})\right)  \\
\nonumber
&& \qquad + \left(\mathbf{U} \otimes \Idvar \right) \, \yglm^{[n-1]}- R_n,\quad  \\ 
&=& h\, \mathbf{S}^{-1}\,\left(\mathbf{A}\E \otimes \Idvar \right) \, \Delta\fex(Y^{[n]}) + h\,\mathbf{S}^{-1}\,\left( \mathbf{A}\I \otimes \Idvar \right) \, \Delta\mathbf{r}\I(Y^{[n]}) \\
\nonumber
&& \qquad + \mathbf{S}^{-1}\,\left(\mathbf{U} \otimes \Idvar \right) \, \Delta\yglm^{[n-1]} - \mathbf{S}^{-1}\,R_n, \qquad \\
R_n &\sim& \mathcal{O}(h^{q\E+1}).
\end{eqnarray*}
Let $h^\ast$ be an upper bound for the step size with the property that
\[
h^\ast < \frac{1}{C_{\mathbf{S}}\,(\Vert \mathbf{A}\E \Vert \, \mathrm{L}\E + \Vert \mathbf{A}\I \Vert \, \mathrm{L}\I)}.
\]
For any $h \le h^\ast$ it holds that:
\begin{eqnarray*}
\Vert \Delta Y^{[n]} \Vert
&\le& h^\ast\,C_{\mathbf{S}}\,(\Vert \mathbf{A}\E \Vert \, \mathrm{L}\E + \Vert \mathbf{A}\I \Vert \, \mathrm{L}\I)\, \Vert \Delta Y^{[n]}\Vert  + C_{\mathbf{S}}\,\Vert \mathbf{U} \Vert \, \Vert \Delta\yglm^{[n-1]} \Vert + C_{\mathbf{S}}\,\Vert R_n\Vert\\ 
&\le&  \alpha \,\Vert \mathbf{U} \Vert \Vert \Delta\yglm^{[n-1]} \Vert + \alpha\,\Vert R_n\Vert, \\
\alpha &=& \left(1 - h^\ast\,C_{\mathbf{S}}\,(\Vert \mathbf{A}\E \Vert \, \mathrm{L}\E + \Vert \mathbf{A}\I \Vert \, \mathrm{L}\I) \right)^{-1} \,C_{\mathbf{S}}.
\end{eqnarray*}
For the final solution the recurrence of global errors is:
\begin{eqnarray*}
\Delta\yglm^{[n]} &=&   \mathbf{M}\I (h\mathbf{J}) \, \Delta\yglm^{[n-1]} + h\, \left(\mathbf{B}\E \otimes \Idvar \right) \, \Delta\fex(Y^{[n]}) + h\,\left( \mathbf{B}\I \otimes \Idvar \right) \, \Delta\mathbf{r}\I(Y^{[n]}) \\
\nonumber
&&\quad + h\,\left( \mathbf{B}\I \otimes h\mathbf{J} \right) \,\mathbf{S}^{-1}\,
\Big( \big(\mathbf{A}\E \otimes \Idvar \big) \, \Delta\fex(Y^{[n]}) \\
&&\qquad + \big( \mathbf{A}\I \otimes \Idvar \big) \, \Delta\mathbf{r}\I(Y^{[n]}) - \mathbf{S}^{-1}\,R_n\Big) - r_n, \\
r_n &\sim& \mathcal{O}(h^{p+1}).
\end{eqnarray*}
Let $\mathbf{U}^{\#}$ be the right pseudo-inverse of $\mathbf{U}$ such that $\mathbf{U}\mathbf{U}^{\#}=\Idstage$, and use the identity
\[
\left( \mathbf{B}\I \otimes h\mathbf{J} \right) \,\mathbf{S}^{-1} = \left( \mathbf{M}\I (h\mathbf{J}) - \mathbf{V}\I \right)\,\mathbf{U}^{\#},
\]
to rewrite the external stage errors as:
\begin{eqnarray*}
\Delta\yglm^{[n]} =   \mathbf{M}\I (h\mathbf{J}) \, \Delta\yglm^{[n-1]} + h\, \left(\mathbf{C}\E \otimes \Idvar\right)\, \Delta\fex(Y^{[n]}) + h\,\left(\mathbf{C}\I \otimes \Idvar\right) \, \Delta\mathbf{r}\I(Y^{[n]}) + \boldsymbol{\rho}_n, 
\end{eqnarray*}
where
\begin{eqnarray*}
\mathbf{C}\E &=&  \mathbf{B}\E + \big( \mathbf{M}\I (h\mathbf{J}) - \mathbf{V}\I \big)\,\mathbf{U}^{\#}\,\mathbf{A}\E, \\
\mathbf{C}\I &=& \mathbf{B}\I + \big( \mathbf{M}\I (h\mathbf{J}) - \mathbf{V}\I \big)\,\mathbf{U}^{\#}\,\mathbf{A}\I,  \\
\boldsymbol{\rho}_n &=& - h\,\big( \mathbf{M}\I (h\mathbf{J}) - \mathbf{V}\I \big)\,\mathbf{U}^{\#}\,\mathbf{S}^{-1}\,R_n - r_n, \\
\Vert \boldsymbol{\rho}_n \Vert &\le& h\,C\,\Vert R_n \Vert + \Vert r_n \Vert \sim \mathcal{O}(h^{\min (p+1,q+2)}) =  \mathcal{O}(h^{p+1}).
\end{eqnarray*}
Solving the recurrence leads to the following global errors:
\begin{eqnarray*}
\Delta\yglm^{[n]}&=&   \left( \mathbf{M}\I (h\mathbf{J}) \right)^n \, \Delta\yglm^{[0]} \\
&& +  \sum_{\ell=1}^{n} \left( \mathbf{M}\I (h\mathbf{J})\right)^{n-\ell} \cdot  \left( h\, \mathbf{C}\E\, \Delta\fex(Y^{[\ell]}) + h\,\mathbf{C}\I \, \Delta\mathbf{r}\I(Y^{[\ell]})  + \boldsymbol{\rho}_\ell \right).
\end{eqnarray*}
Taking norms we have, for $h \le h^\ast$:
\begin{eqnarray*}
\Vert \Delta\yglm^{[n]} \Vert &\le&   C_{\mathbf{M}} \, \Delta\yglm^{[0]} + \underbrace{ C_{\mathbf{M}} \,
 \left( \Vert \mathbf{C}\E \Vert\, \mathrm{L}\E  +\Vert \mathbf{C}\I \Vert\, \mathrm{L}\I \right)}_{\beta}\,\sum_{\ell=1}^{n} h\, \Vert \Delta Y^{[\ell]}\Vert,\\
&& \quad + C_{\mathbf{M}} \,\sum_{\ell=1}^{n}  \boldsymbol{\rho}_\ell \\
&\le& C_{\mathbf{M}} \, \Delta\yglm^{[0]} + \underbrace{\beta\,\alpha\,\Vert \mathbf{U}\Vert}_{\eqqcolon \zeta}\,\sum_{\ell=1}^{n} h\, \Vert \Delta\yglm^{[\ell-1]} \Vert + \sum_{\ell=1}^{n} \underbrace{\left( C_{\mathbf{M}} \, \boldsymbol{\rho}_\ell + \alpha\,h\,\norm{ R_\ell } \right)}_{\eqqcolon \delta^{[\ell]} \sim \mathcal{O}\big( h^{p+1} \big)}.
\end{eqnarray*}
By standard techniques we compute the sequence of upper bounds $\eta^{[n]}$ defined by
\begin{eqnarray*}
\eta^{[n]} &=&  C_{\mathbf{M}} \, \eta^{[0]} +   \zeta \sum_{\ell=1}^{n-1} \left( h\,\eta^{[\ell]} + \delta^{[\ell]} \right),\quad
\delta^{[\ell]} \sim \mathcal{O}\big( h^{p+1} \big).
\end{eqnarray*}
From the recurrences we have that $\Vert \Delta\yglm^{[n]} \Vert \le \eta^{[n]}$ for all $n$. The $\eta^{[n]}$ recurrence gives:
\[
\eta^{[n]} = (1 + h\,\zeta)\, \eta^{[n-1]} + \zeta\, \delta^{[\ell]} \quad \Rightarrow \quad \eta^{[n]} \le e^{(t_F-t_0)\,\zeta} \,
\left( \eta^{[0]} + \zeta \sum_{\ell=0}^{n-1} \delta^{[\ell]} \right),
\]
for $nh=t_F-t_0$.
%
By standard arguments this stability relation,  the local truncation errors of size $\delta^{[\ell]} \sim \mathcal{O}(h^{p+1})$, and the initial error $\eta^{[0]} \sim \mathcal{O}(h^{p})$, imply a global error of size $\mathcal{O}(h^{p})$.

\end{proof}

\section{Singular perturbation and index-1 differential-algebraic problems}
\label{sec:problems}

\subsection{The index-1 differential algebraic problem}

Consider the index-1 differential algebraic equation (DAE) \cite{Hairer_DAE_book,Hairer_book_II,Tikhonov_1985_ODE}:
\begin{equation}
\label{eqn:DAE-index1}
\begin{cases}
\xx' = \fun(\xx,\zz), \\
0 = \gun(\xx,\zz),
\end{cases}
\end{equation}
where $\fun$, $\gun$ are smooth functions and the sub-Jacobian $\gun_\zz$ is invertible in a neighborhood of the solution. 
The initial values $[\xx_0,\zz_0]$ are consistent if $\gun(\xx_0,\zz_0)=0$.
By the implicit function theorem the algebraic equation can be locally solved uniquely to express $\zz$ as a function of $\xx$:
\[
\zz = \mathcal{G}(\xx).
\]
Replacing this in the differential equation \eqref{eqn:DAE-index1} leads to the following reduced ODE:
\begin{equation}
\label{eqn:reduced-ode}
\xx' = \fun\big(\xx,\mathcal{G}(\xx)\big) \eqqcolon \fun^{\textsc{red}}\big(\xx\big).
\end{equation}

\subsection{The singular perturbation problem}

Consider the singular perturbation problem \cite{Hairer_DAE_book,Hairer_book_II,Tikhonov_1985_ODE}
\begin{equation}
\label{eqn:SPP}
\begin{cases}
\xx' = \fun(\xx,\zz), \\
\zz' = \varepsilon^{-1}\,\gun(\xx,\zz), 
\end{cases}
\end{equation}
where $\varepsilon \ll 1$. The Jacobian $\gun_\zz$ is assumed to be invertible and with a negative logarithmic norm 
\begin{equation}
\label{eqn:SPP-logarithmic-norm}
\mu\left( \gun_\zz(\xx,\zz) \right)  = \lim_{h \to 0+} \frac{\Vert \Id + h\, \gun_\zz(\xx,\zz) \Vert -1 }{h} = \sup_{\yy \ne 0} \frac{Re \, \langle \gun_\zz(\xx,\zz)\,\yy, \yy \rangle }{\langle \yy, \yy \rangle}  \le -1
\end{equation}
in an $\varepsilon$-independent neighborhood of the solution. (Any other negative bound can be scaled to $-1$). Consequently, in the limit  $\varepsilon \to 0$ the system \eqref{eqn:SPP} becomes an index-1 DAE \eqref{eqn:DAE-index1}.

\begin{remark}[Additively partitioned form]
The system \eqref{eqn:SPP} can be written as an additively partitioned system:
\begin{equation}
\label{eqn:SPP-additive}
\yy \coloneqq \begin{bmatrix} \xx \\ \zz \end{bmatrix}, \qquad
\underbrace{\begin{bmatrix} \xx \\ \zz \end{bmatrix}'}_{\yy'} = 
\underbrace{\begin{bmatrix} \fun(\xx,\zz) \\ 0 \end{bmatrix}}_{\fex(\yy)} + 
\underbrace{\begin{bmatrix} 0 \\ \varepsilon^{-1}\,\gun(\xx,\zz) \end{bmatrix}}_{\fim(\yy)}. 
\end{equation}
\end{remark}

For arbitrary initial values the solution of \eqref{eqn:SPP} undergoes a short transient over an $\mathcal{O}(\varepsilon)$ long time interval; this transient dies out, and afterwards the solution of  \eqref{eqn:SPP} is smooth and stiff \cite{Tikhonov_1985_ODE}.
We assume that the initial conditions are along the smooth solution (i.e., we consider the solution after the transient has passed). The solutions of the singular perturbation problem \eqref{eqn:SPP} can then be expanded in a series of powers of $\varepsilon$ \cite{Hairer_DAE_book}:
\begin{equation}
\label{eqn:SPP-epsilon-expansion}
\xx(t) = \sum_{k \ge 0} \xx^{k}(t)\,\varepsilon^k, \quad
\zz(t) = \sum_{k \ge 0} \zz^{k}(t)\,\varepsilon^k, 
\end{equation}
where $\xx^{k}(t)$ and $\zz^{k}(t)$ are smooth functions that do not depend on $\varepsilon$.

Inserting \eqref{eqn:SPP-epsilon-expansion} into \eqref{eqn:SPP}, and equating powers of $\varepsilon$, shows that these functions are solutions of differential-algebraic equations of increasing indices, as follows \cite{Hairer_DAE_book,Hairer_book_II}. 
Specifically, we insert \eqref{eqn:SPP-epsilon-expansion}  into \eqref{eqn:SPP}  and expand the function in Taylor series about the $\mathcal{O}(1)$ terms to obtain:
\begin{subequations}
\label{eqn:spp-eps-expansion}
\begin{eqnarray}
\sum_{k \ge 0} (\xx^{k})'\,\varepsilon^k &=& \fun(\xx^{0},\zz^{0}) + 
\sum_{\alpha+\beta\ge 1} \frac{\partial^{\alpha+\beta}\fun}{\partial\xx^\alpha\partial\zz^\beta}(\xx^{0},\zz^{0}) \cdot \big( W_{\alpha,\beta} ) \\
\sum_{k \ge 0} (\zz^{k})'\,\varepsilon^{k+1} &=& \gun(\xx^{0},\zz^{0}) + 
\sum_{\alpha+\beta\ge 1} \frac{\partial^{\alpha+\beta}\gun}{\partial\xx^\alpha\partial\zz^\beta}(\xx^{0},\zz^{0}) \cdot \big( W_{\alpha,\beta} ),
\end{eqnarray}
\end{subequations}
where the function derivatives are applied to the following multinomial arguments:
\begin{equation*}
\begin{split}
\big( W_{\alpha,\beta} \big) &\coloneqq  \Big( \underbrace{\dots, \sum_{m \ge 1} \xx^{m}\,\varepsilon^m, \dots}_{\alpha\textnormal{ times}}~, ~ \underbrace{\dots,\sum_{\ell \ge 1} \zz^{\ell}\,\varepsilon^\ell, \dots}_{\beta\textnormal{ times}} \Big) = \sum_{k \ge 1} \big( W_{\alpha,\beta}^{k} \big)\, \varepsilon^{k}, \\
\big( W_{\alpha,\beta}^{k} \big) &= \sum_{\sum_{i = 1}^\alpha m_i + \sum_{j = 1}^\beta \ell_j = k} 
\Big( \xx^{m_1},\dots, \xx^{m_\alpha}, \zz^{\ell_1} , \dots, \zz^{\ell_\beta} \Big).
\end{split}
\end{equation*}

%
Equating the coefficients of $\varepsilon^0$ in \eqref{eqn:spp-eps-expansion} gives:
\begin{subequations}
\label{eqn:dae-eps-0}
\begin{eqnarray}
\label{eqn:dae-eps-0-x}
(\xx^0)' &= \fun(\xx^0,\zz^0), \\
\label{eqn:dae-eps-0-z}
0 &= \gun(\xx^0,\zz^0). 
\end{eqnarray}
\end{subequations}
Consequently, the $\mathcal{O}(1)$ terms $(\xx^0,\zz^0)$ are the solution of the index-1 DAE system \eqref{eqn:DAE-index1}.
Differentiating  \eqref{eqn:dae-eps-0-z} gives an ODE for the algebraic variable $\zz^0(t)$:
\begin{equation}
(\zz^0)' = \mathcal{G}_\xx(\xx^0)\,(\xx^0)' = - (\gun_\zz^{-1}\,\gun_\xx\,\fun)(\xx^0,\zz^0).
\end{equation}
Equating the coefficients of $\varepsilon^{\nu}$, $\nu \ge 1$, in \eqref{eqn:spp-eps-expansion} gives:
\begin{subequations}
\label{eqn:spp-eps-expansion-k}
\begin{eqnarray}
\label{eqn:spp-eps-expansion-k-x}
(\xx^{\nu})'(t) &=& 
\sum_{1 \le \alpha+\beta \le \nu} \frac{\partial^{\alpha+\beta}\fun}{\partial\xx^\alpha\partial\zz^\beta}(\xx^{0},\zz^{0}) \cdot \big( W_{\alpha,\beta}^{\nu}  ) \\
\nonumber
&=& \fun_\xx(\xx^0,\zz^0)\,\xx^\nu + \fun_\zz(\xx^0,\zz^0)\,\zz^\nu + \sum_{2 \le \alpha+\beta \le \nu} \frac{\partial^{\alpha+\beta}\fun}{\partial\xx^\alpha\partial\zz^\beta}(\xx^{0},\zz^{0}) \cdot \big( W_{\alpha,\beta}^{\nu}  ), \\
\label{eqn:spp-eps-expansion-k-z}
(\zz^{\nu-1})'(t) &=& \
\sum_{1 \le \alpha+\beta \le \nu} \frac{\partial^{\alpha+\beta}\gun}{\partial\xx^\alpha\partial\zz^\beta}(\xx^{0},\zz^{0}) \cdot \big( W_{\alpha,\beta}^{\nu} ) \\
\nonumber
&=& \gun_\xx(\xx^0,\zz^0)\,\xx^\nu + \gun_\zz(\xx^0,\zz^0)\,\zz^\nu + \sum_{2 \le \alpha+\beta \le \nu} \frac{\partial^{\alpha+\beta}\gun}{\partial\xx^\alpha\partial\zz^\beta}(\xx^{0},\zz^{0}) \cdot \big( W_{\alpha,\beta}^{\nu} )
\end{eqnarray}
\end{subequations}
Solving for the algebraic variable $\zz^{\nu}(t)$ gives:
\begin{subequations}
\label{eqn:expansion-eps-nu}
\begin{eqnarray}
\label{eqn:expansion-eps-nu-X}
(\xx^\nu)' &=& (\fun_\xx - \fun_\zz\,\gun_\zz^{-1}\,\gun_\xx)\,\xx^{\nu} +(\fun_\zz\,\gun_\zz^{-1})\,(\zz^{\nu-1})' + \sum_{2 \le \alpha+\beta \le \nu} \frac{\partial^{\alpha+\beta}\fun}{\partial\xx^\alpha\partial\zz^\beta} \cdot \big( W_{\alpha,\beta}^{\nu}  ),\qquad\quad \\
\label{eqn:expansion-eps-nu-Z}
\zz^\nu &=& - (\gun_\zz^{-1}\,\gun_\xx)\,\xx^{\nu} + \gun_\zz^{-1}\,(\zz^{\nu-1})' - \gun_\zz^{-1}\, \sum_{2 \le \alpha+\beta \le \nu} \frac{\partial^{\alpha+\beta}\gun}{\partial\xx^\alpha\partial\zz^\beta} \cdot \big( W_{\alpha,\beta}^{\nu} ),
\end{eqnarray}
\end{subequations}
where all function derivatives are evaluated at $(\xx^0,\zz^0)$.

Equating \eqref{eqn:spp-eps-expansion-k} for $\nu=1$ leads to
\begin{subequations}
\label{eqn:dae-eps-1}
\begin{eqnarray}
\label{eqn:dae-eps-1-X}
(\xx^1)' &= \fun_\xx(\xx^0,\zz^0)\,\xx^1 + \fun_\zz(\xx^0,\zz^0)\,\zz^1, \\
\label{eqn:dae-eps-1-Z}
(\zz^0)' &= \gun_\xx(\xx^0,\zz^0)\,\xx^1 + \gun_\zz(\xx^0,\zz^0)\,\zz^1.
\end{eqnarray}
\end{subequations}
Note that this is a linear index-2 DAE in $(\xx^1(t),\zz^1(t))$. Solving for the algebraic variable $\zz^1(t)$ gives:
\begin{subequations}
\label{eqn:expansion-eps1}
\begin{eqnarray}
\label{eqn:expansion-eps1-X}
(\xx^1)' &=& (\fun_\xx - \fun_\zz\,\gun_\zz^{-1}\,\gun_\xx)\,\xx^1 +(\fun_\zz\,\gun_\zz^{-1})\,(\zz^0)', \\
\label{eqn:expansion-eps1-Z}
\zz^1 &=& - (\gun_\zz^{-1}\,\gun_\xx)\,\xx^1 + \gun_\zz^{-1}\,(\zz^0)',
\end{eqnarray}
\end{subequations}
where all Jacobians are evaluated at $(\xx^0,\zz^0)$.

\section{Existence and uniqueness of the IMEX GLM numerical solution on singular perturbation problem}
\label{sec:SPP-solution-existence}
%
Schneider \cite{Schneider_1993_GLM-SPP} performed a detailed singular perturbation analysis of implicit (non-partitioned) GLMs. 
Other authors also studied implicit (non-partitioned) GLMs applied to differential-algebraic problems. 
Chartier \cite{Chartier_1993_GLM-DAE} studied the convergence of GLMs for DAEs of index-1, and of stiffly accurate GLMs for DAEs of index-2. Butcher and Chartier \cite{Butcher_1995_GLM-DAE} developed parallel implicit GLMs and applied them to stiff ordinary differential and differential algebraic equations of index two and three.
Schulz \cite{Schulz_2003_GLM-DAE} derived order conditions and starting methods for stiffly accurate GLMs applied to {\it linear} DAEs with properly stated leading terms of index 1 and 2.
Voigtmann \cite{Voigtmann_2005_GLM-DAE} studied stiffly accurate GLMs to solve nonlinear index-2 DAEs with a special structure appearing in the 
modeling of electrical circuits.

We solve the singular perturbation problem \eqref{eqn:SPP-additive} with the IMEX-GLM method \eqref{eqn:component-imex-glm} in component partitioned form, where the $\xx$ variables are non-stiff, and the $\zz$ variables are stiff: 
\begin{subequations}
\label{eqn:IMEX_SPP}
\begin{eqnarray}
\label{eqn:IMEX_SPP-intX}
X^{[n]} &=& h\, \left(\mathbf{A}\E \otimes \Idvar \right) \, \fun\big(X^{[n]},Z^{[n]}\big)  + \left(\mathbf{U}\E \otimes \Idvar \right) \, \xglm^{[n-1]},\qquad \\ 
\label{eqn:IMEX_SPP-intZ}
Z^{[n]} &=& h\,\left( \mathbf{A}\I \otimes \Idvar \right) \, \varepsilon^{-1}\,\gun\big(X^{[n]},Z^{[n]}\big)  + \left(\mathbf{U}\I \otimes \Idvar \right) \, \zglm^{[n-1]},\qquad \\ 
\label{eqn:IMEX_SPP-extX}
\xglm^{[n]} &=& h\, \left(\mathbf{B}\E \otimes \Idvar \right) \, \fun\big(X^{[n]},Z^{[n]}\big)  + \left(\mathbf{V}\E \otimes \Idvar \right) \, \xglm^{[n-1]}, \\
\label{eqn:IMEX_SPP-extZ}
\zglm^{[n]} &=& h\,\left( \mathbf{B}\I \otimes \Idvar \right)  \, \varepsilon^{-1}\,\gun\big(X^{[n]},Z^{[n]}\big)  + \left(\mathbf{V}\I \otimes \Idvar \right) \, \zglm^{[n-1]}.\qquad
\end{eqnarray}
\end{subequations}

For small step sizes $h < \varepsilon$ standard arguments show that the nonlinear system of equations \eqref{eqn:IMEX_SPP-intX}--\eqref{eqn:IMEX_SPP-intZ} has a unique solution. We show that the system also has a unique solution for the very stiff case where $h \gg \varepsilon$.

\begin{theorem}[Existence of IMEX GLM solution on SPP problem]
\label{thm:solution-existence}
Assume that \eqref{eqn:SPP-logarithmic-norm} $\mu\left( \gun_\zz(\xx,\zz) \right) \le -1$ holds, and that the eigenvalues of the $\mathbf{A}\I$ matrix have positive real part. Assume that the external stages at time $t_{n-1}$ are order one accurate:
\begin{equation}
\label{eqn:initial-conditions-assumption}
\xglm^{[n-1]} = \mathbf{w}_0\, \xx(t_{n-1}) + \mathcal{O}(h), \qquad 
\zglm^{[n-1]} = \mathbf{w}_0\, \zz(t_{n-1}) + \mathcal{O}(h).
\end{equation}
Then  the nonlinear system of equations \eqref{eqn:IMEX_SPP-intX}--\eqref{eqn:IMEX_SPP-intZ} has a unique solution for any $h \le h_0$, with $h_0$ sufficiently small but independent of $\varepsilon$. The solution satisfies:
\[
X^{[n]} = \yy(t_{n-1}) + \mathcal{O}(h), \qquad Z^{[n]} = \zz(t_{n-1}) + \mathcal{O}(h),
\]
\end{theorem}

\begin{proof}
The reasoning follows closely the proofs of existence of a numerical solution for implicit Runge-Kutta methods \cite{Hairer_1988_RK-SPP}, \cite[Theorem VI.3.5]{Hairer_book_II}, and for implicit GLMs \cite{Schneider_1993_GLM-SPP}. Let
\begin{equation}
\label{eqn:initial-conditions}
\begin{split}
\eta &\coloneqq \left(\mathbf{U}\E \otimes \Idvar \right) \, \xglm^{[n-1]} 
= \one_s \otimes \xx(t_{n-1}) + \mathcal{O}(h), \\
\zeta &\coloneqq \left(\mathbf{U}\I \otimes \Idvar \right) \, \zglm^{[n-1]} = \one_s \otimes \zz(t_{n-1}) + \mathcal{O}(h), 
\end{split}
\end{equation}
where we use the assumption \eqref{eqn:initial-conditions-assumption} together with the preconsistency condition \eqref{eqn:GLM-preconsistency}. Consequently $\gun( \eta,\zeta ) = \mathcal{O}(h)$.

The nonlinear system \eqref{eqn:IMEX_SPP-intX}--\eqref{eqn:IMEX_SPP-intZ} can be written as:
\begin{equation}
\label{eqn:existence-nonlinear-system}
\begin{split}
X^{[n]} - \eta &= h\, \left(\mathbf{A}\E \otimes \Idvar \right) \, \fun\big(X^{[n]},Z^{[n]}\big),\qquad \\ 
(\varepsilon/h)\, \left( Z^{[n]} - \zeta \right) &= \left( \mathbf{A}\I \otimes \Idvar \right) \, \gun\big(X^{[n]},Z^{[n]}\big),
\end{split}
\end{equation}
and has a Jacobian of the form:
\[
\begin{bmatrix}
\Id + \mathcal{O}(h) & \mathcal{O}(h) \\
\mathcal{O}(1) & (\varepsilon/h)\,\Id - \mathbf{A}\I \otimes \gun\big(X^{[n]},Z^{[n]}\big)
\end{bmatrix}.
\]
The proof follows exactly as in \cite[Theorem VI.3.5]{Hairer_book_II} by applying Newton iterations from the initial point $(\eta,\zeta)$, showing that the Jacobian has a uniformly bounded inverse, and applying the Newton-Kantorovich theorem.
\end{proof}

\section{Convergence analysis for index-1 differential-algebraic problems}
\label{sec:SPP-index-1}
%
Schneider \cite{Schneider_1993_GLM-SPP} performed a detailed singular perturbation analysis of implicit (non-partitioned) GLMs. He provided convergence results for index-1 DAEs, and for index-2 DAEs for a special class of GLMs. Here, and in the next section, we extend his analysis to IMEX GLMs.

\subsection{Accumulation of errors}

We start our analysis with the following Lemma.

\begin{lemma}[Accumulation of errors]
\label{lem:error-convergence}
Consider the iteration
\begin{equation}
\label{eqn:skeleton-iteration}
\zeta^{[n]} = \mathbf{M} \, \zeta^{[n-1]} + \delta^{[n]} = \sum_{k=0}^n \mathbf{M}^{n-k} \, \delta^{[k]}, \quad \zeta^{[0]} =  \delta^{[0]}, \quad \delta^{[k]} \sim \mathcal{O}\bigl( h^{\nu} \bigr).
\end{equation}
The  iteration matrix is power bounded and has a Jordan form:
\[
\sup_n\,\Vert \mathbf{M}^n \Vert = C < \infty, \qquad
\mathbf{X}\, \mathbf{M} \, \mathbf{X}^{-1} = \begin{bmatrix} \mathbf{D} & \Zero & \Zero & \\ \Zero & \boldsymbol{\Lambda} & \Zero \\ \Zero & \Zero & \mathbf{J} \end{bmatrix},
\]
where:
\begin{enumerate}
\item $\boldsymbol{D}$ is a diagonal matrix whose diagonal entries are the regular eigenvalues of $\mathbf{M}$ with absolute value of one and of the form $e^{i \pi /L}$ with $L$ an integer; the case $L=1$ corresponds to an eigenvalue of $-1$;
\item $\boldsymbol{\Lambda}$ is a diagonal matrix whose diagonal entries are the regular eigenvalues of $\mathbf{M}$ with absolute value of one but not equal to $e^{i \pi /L}$; 
\item $\mathbf{J} = \operatorname{blkdiag} \mathbf{J}_k$ is block diagonal, with each Jordan block $\mathbf{J}_k$ corresponding to an eigenvalue of $\mathbf{M}$ with absolute value strictly less than one. 
 \end{enumerate}
Apply the Jordan transform to the iteration \eqref{eqn:skeleton-iteration}, and denote 
 \[
 \begin{bmatrix} \zeta^{[n]}_1 \\ \zeta^{[n]}_2  \\ \zeta^{[n]}_3 \end{bmatrix} = \mathbf{X}\, \zeta^{[n]}, \quad
 \begin{bmatrix} \delta^{[n]}_1 \\ \delta^{[n]}_2  \\ \delta^{[n]}_3 \end{bmatrix} = \mathbf{X}\, \delta^{[n]}.
 \]
 We have the following convergence result:
 \begin{subequations}
\begin{eqnarray}
\label{eqn:lemma-case-1}
 \zeta_1^{[n]} &\sim& \begin{cases} \mathcal{O}\bigl( h^{\nu-1} \bigr) & \textnormal{in general}, \\
 \label{eqn:lemma-case-2}
\mathcal{O}\bigl( h^{\nu} \bigr) & \textnormal{if } \delta^{[i]}_1-\delta^{[i-1]}_1 \sim \mathcal{O}\bigl( h^{\nu+1} \bigr), \end{cases}  \\
 \zeta_2^{[n]} &\sim& \mathcal{O}\bigl( h^{\nu-1} \bigr), \\
\label{eqn:lemma-case-3}
 \zeta_3^{[n]} &\sim& \mathcal{O}\bigl( h^{\nu} \bigr).
 \end{eqnarray}
 \end{subequations}
\end{lemma}
\begin{proof}

We consider the iterations \eqref{eqn:skeleton-iteration} for individual components. Since $\mathbf{D}$ is diagonal, the analysis of the iterations for the first component $\zeta_1^{[n]}$ can be done component-wise. Without loss of generality we consider a single eigenvalue $e^{i \pi/L}$. We have:
\begin{equation*}
\begin{split}
\zeta_1^{[n]} &= \sum_{k=0}^n e^{i \pi (n-k)/L} \, \delta_1^{[k]} = e^{i \pi \,n/L} \,\sum_{m\,\textnormal{odd}} \sum_{\ell=0}^{L-1}   e^{-i \pi \ell/L} \, (- \delta^{[Lm+\ell]} + \delta^{[L(m-1)+\ell]}), \\
\Vert \zeta_1^{[n]} \Vert &\le 
\begin{cases} \sum_{k=0}^n \Vert \delta_1^{[k]} \Vert \sim \mathcal{O}\bigl( h^{\nu-1} \bigr) & \textnormal{in general}, \\ 
\sum_{m\,\textnormal{odd}} \sum_{\ell=0}^{L-1}  \Vert \delta_1^{[Lm+\ell]}-\delta_1^{[Lm+\ell-L]} \Vert \sim \mathcal{O}\bigl( h^{\nu} \bigr) & \textnormal{if } \delta^{[k]}_1-\delta^{[k-L]}_1 \sim \mathcal{O}\bigl( h^{\nu+1} \bigr).
\end{cases} \\
\end{split}
\end{equation*}
Of special interest is the case $L=1$ when the eigenvalues are equal to $-1$.

For the second components we use the fact that  $\boldsymbol{\Lambda}$ is a diagonal matrix with diagonal entries of absolute value equal to one. We have:
\begin{equation*}
\begin{split}
\Vert \zeta_2^{[n]} \Vert &\le \sum_{k=0}^n \Vert \boldsymbol{\Lambda}^{n-k} \Vert \,\Vert \delta^{[k]} \Vert
= \sum_{k=0}^n \Vert \delta^{[k]} \Vert \sim \mathcal{O}\bigl( h^{\nu-1} \bigr), \\
\end{split}
\end{equation*}

To analyze the third components we first note that  $\Vert \mathbf{J}_k^n \Vert \le C$ for any power $n \ge 0$, and that: 
\[
\exists\, n_0,\rho \textnormal{ with } 0 < \rho < 1 \textnormal{ such that: } \Vert \mathbf{J}_k^n \Vert \le \rho < 1~~\forall\, n \ge n_0.
\]
We have:
\begin{equation*}
\begin{split}
\Vert \zeta_3^{[n]} \Vert &\le \sum_{k=0}^n \Vert \mathbf{J}^{n-k} \Vert \,\Vert \delta^{[k]} \Vert
\le \sum_{k=0}^{n_0-1} \Vert \mathbf{J}^{n-k} \Vert\,\Vert \delta^{[k]} \Vert + \sum_{k=n_0}^n 
\rho^{n-k} \Vert \delta^{[k]} \Vert \\
&\le \left( n_0\,C + \frac{1}{1-\rho}\right)\,\max_k\, \Vert \delta^{[k]} \Vert \sim \mathcal{O}\bigl( h^{\nu} \bigr).
\end{split}
\end{equation*}
\end{proof}

\subsection{The index-1 DAE solution}
\label{sec:SPP-index-1}
Consider the IMEX GLM  solution for the singular perturbation problem \eqref{eqn:IMEX_SPP}, and take the limit $\varepsilon \to 0$. The method has an invertible $\mathbf{A}\I$. Taking the limit $\varepsilon \to 0$ in \eqref{eqn:IMEX_SPP-intZ} leads to:
\begin{equation}
\label{eqn:DAE-stages-in-the-limit}
\gun\big(X^{[n]},Z^{[n]}\big) = \Zero.
\end{equation}
From \eqref{eqn:IMEX_SPP-intZ}, using the fact that $\mathbf{A}\I$ is invertible, express stage values of the stiff function as
\[
h\, \varepsilon^{-1}\,\gun\big(X^{[n]},Z^{[n]}\big) = \left( \mathbf{A}\I\,^{-1} \otimes \Idvar \right)\, \left(Z^{[n]} - \left(\mathbf{U}\I \otimes \Idvar \right) \, \zglm^{[n-1]} \right),  
\]
and substitute in \eqref{eqn:IMEX_SPP-extZ} to obtain:
\begin{equation}
\label{eqn:stiff-stages}
\begin{split}
\zglm^{[n]} &= \left( \mathbf{B}\I\,\mathbf{A}\I\,^{-1} \otimes \Idvar \right)  \,  Z^{[n]}  + \Bigl( \underbrace{\bigl( \mathbf{V}\I - \mathbf{B}\I\,\mathbf{A}\I\,^{-1}\,\mathbf{U}\I  \bigr) }_{\mathbf{M}\I(\infty)} \otimes \Idvar \Bigr)\, \zglm^{[n-1]}.
\end{split}
\end{equation}
In the limit $\varepsilon \to 0$ the IMEX-GLM provides the following numerical solution of the index-1 DAE problem \eqref{eqn:DAE-index1}: 
\begin{subequations}
\label{eqn:IMEX_DAE1}
\begin{eqnarray}
\label{eqn:IMEX_DAE1-intX}
X^{[n]} &=& h\, \left(\mathbf{A}\E \otimes \Idvar \right) \, \fun\big(X^{[n]},Z^{[n]}\big)  + \left(\mathbf{U}\E \otimes \Idvar \right) \, \xglm^{[n-1]},\qquad \\ 
\label{eqn:IMEX_DAE1-intZ}
0 &=& \gun\big(X^{[n]},Z^{[n]}\big),\qquad \\ 
\label{eqn:IMEX_DAE1-extX}
\xglm^{[n]} &=& h\, \left(\mathbf{B}\E \otimes \Idvar \right) \, \fun\big(X^{[n]},Z^{[n]}\big)  + \left(\mathbf{V}\E \otimes \Idvar \right) \, \xglm^{[n-1]}, \\
\label{eqn:IMEX_DAE1-extZ}
\zglm^{[n]} &=& \left( \mathbf{B}\I\,\mathbf{A}\I\,^{-1} \otimes \Idvar \right)  \,  Z^{[n]}  + \Bigl( \mathbf{M}\I(\infty) \otimes \Idvar \Bigr)\, \zglm^{[n-1]}.\qquad
\end{eqnarray}
\end{subequations}
We denote the global errors of the numerical solution \eqref{eqn:IMEX_DAE1} by:
\begin{equation}
\label{eqn:global-error-notation}
\begin{split}
\Delta\xglm^{[n]} &\coloneqq \xglm^{[n]} - \left(\mathbf{W}\E\otimes \Idvar \right)\,\Nordsieck_p(h,\xx,t_n), \\
\Delta\zglm^{[n]} &\coloneqq \zglm^{[n]} - \left(\mathbf{W}\I\otimes \Idvar \right)\,\Nordsieck_p(h,\zz,t_n),
\end{split}
\end{equation}
and the global stage errors by
\begin{equation}
\label{eqn:stage-error-notation}
\begin{split}
\Delta X^{[n]} &\coloneqq X^{[n]} - \xx(t_{n-1}+\mathbf{c} h), \\
\Delta Z^{[n]} &\coloneqq Z^{[n]} - \zz(t_{n-1}+\mathbf{c} h).
\end{split}
\end{equation}

\begin{theorem}[Order of IMEX-GLM on explicit index-1 DAE]
\label{thm:index-1-order}
If the implicit component has a coefficient matrix $\mathbf{A}\I$ with (strictly) positive eigenvalues then the non-stiff external stages converge with order $p$, the stiff external stages converge with order $\nu$,
\begin{subequations}
\label{eqn:index-1-order}
\begin{align}
\Delta X^{[n]} & \sim \mathcal{O}(h^{\min(p,q\E+1)}), &
\Delta\xglm^{[n]} &\sim \mathcal{O}(h^p), \\
\Delta Z^{[n]} &\sim \mathcal{O}(h^{\min(p,q\E+1)}), &
\Delta\zglm^{[n]} &\sim \mathcal{O}(h^\nu), 
\end{align}
\end{subequations}
where, with $q = \min\{ q\E,q\I \}$, the order of convergence of the stiff variables is:
\begin{enumerate}
\item $\nu = q$ if the implicit stability matrix at infinity is power bounded, $\Vert \bigl(\mathbf{M}\I(\infty)\bigr)^{k} \Vert \le C$ for all $k$; and
\item $\nu = \min(p,q+1)$ if the spectrum of $\mathbf{M}\I(\infty)$ contains eigenvalues of absolute value smaller than one, and the regular eigenvalues of modulus one, if any, have the form $e^{-i \pi/L}$ for $L$ an integer.
\end{enumerate}
\end{theorem}

\begin{remark}
Theorem \ref{thm:index-1-order} is similar to the convergence results of Schneider  \cite{Schneider_1993_GLM-SPP}  and Chartier \cite{Chartier_1993_GLM-DAE} for implicit (non-partitioned) GLMs applied to index-1 DAEs.
\end{remark}

\begin{proof}
Equation \eqref{eqn:IMEX_DAE1-intZ} leads to the condition:
\begin{equation}
\label{eqn:stages-in-the-limit}
\gun\big(X^{[n]},Z^{[n]}\big) = \Zero \quad \Rightarrow \quad Z^{[n]} = \mathcal{G}(X^{[n]}).
\end{equation}
Consequently, equations \eqref{eqn:IMEX_SPP-intX}, \eqref{eqn:stages-in-the-limit}, and \eqref{eqn:IMEX_SPP-extX} represent the explicit component method applied to the reduced system of ODEs \eqref{eqn:reduced-ode}
\begin{eqnarray*}
X^{[n]} &=& h\, \left(\mathbf{A}\E \otimes \Idvar \right) \, \fun^{\textsc{red}}\big(X^{[n]}\big)  + \left(\mathbf{U}\E \otimes \Idvar \right) \, \xglm^{[n-1]},\qquad \\ 
\xglm^{[n]} &=& h\, \left(\mathbf{B}\E \otimes \Idvar \right) \, \fun^{\textsc{red}}\big(X^{[n]}\big)  + \left(\mathbf{V}\E \otimes \Idvar \right) \, \xglm^{[n-1]}.
\end{eqnarray*}
Therefore, the solution $\xglm^{[n]}$ is an approximation of global order $p$, and the stages $X^{[n]}$ are approximations of order $q\E$:
\begin{subequations}
\label{eqn:index1-stage-global-err}
\begin{equation}
\label{eqn:nonstiff-stage-global-err}
\begin{split}
\xglm^{[n]} - \mathbf{W}\E\,\Nordsieck_p(h,\xx,t_n) &\sim \mathcal{O}(h^{p}), \\
X^{[n]} - \xx(t_n + \mathbf{c}\,h) &\sim \mathcal{O}(h^{\min(p,q\E+1)}).
\end{split}
\end{equation}
It follows from \eqref{eqn:stages-in-the-limit} that the stiff stages are also approximations of order $q\E$:
\begin{equation}
\label{eqn:stiff-stage-global-err}
Z^{[n]} = \mathcal{G}\big(\xx(t_{n-1} + \mathbf{c}\,h\big) + \mathcal{O}(h^{q\E+1})) = \zz(t_{n-1} + \mathbf{c}\,h) + \mathcal{O}(h^{\min(p,q\E+1)}).
\end{equation}
\end{subequations}

Consider now the stiff external stage equation \eqref{eqn:IMEX_DAE1-extZ}.
From the preconsistency conditions \eqref{eqn:GLM-preconsistency}
\begin{equation}
\label{eqn:GLMO1}
\mathbf{M}\I(\infty)\,\mathbf{w}_0 
= \mathbf{w}_0 - \mathbf{B}\I\,\mathbf{A}\I\,^{-1}\,\mathbf{1}_{s \times 1}.
\end{equation}
Multiplying the implicit component stage order condition \eqref{eqn:GLM-order_condition-1} by $\mathbf{B}\I\mathbf{A}\I\,^{-1}$ from the left, and subtracting the result from the implicit order condition \eqref{eqn:GLM-order_condition-2} gives:
\begin{equation}
\label{eqn:GLMO2}
 \sum_{\ell=0}^{k}  \frac{\mathbf{w}_{k-\ell}}{\ell!} 
- \frac{\mathbf{B}\I\,\mathbf{A}\I\,^{-1}\,\mathbf{c}^{\times k}}{k!} = \mathbf{M}\I(\infty)\,\mathbf{w}\I_k, \quad
 k=1,\dots,q\I.
\end{equation}
From \eqref{eqn:GLMO1} and \eqref{eqn:GLMO2} for $k=0,\dots,q\I$ we have:
\[
\begin{split}
&  \left( \mathbf{M}\I(\infty)\otimes \Idvar \right) \,\underbrace{\sum_{k =0}^{q\I} \mathbf{w}\I_k \,h^k\, \zz^{(k)}(t_{n-1})}_{\mathbf{W}\I_{:,0:q\I}\, \Nordsieck_{q\I}(h,\zz,t_{n-1})} \\
&= \underbrace{ \sum_{k=0}^{q\I} \sum_{\ell=0}^{k}  \frac{\mathbf{w}\I_{k-\ell}}{\ell!} \,h^k\, \zz^{(k)}(t_{n-1}) }_{\mathbf{W}\I_{:,0:q\I}\, \Nordsieck_{q\I}(h,\zz,t_n)}
  - \mathbf{B}\I\,\mathbf{A}\I\,^{-1}\,\underbrace{ \sum_{k =0}^{q\I} \frac{h^k\,\mathbf{c}^{\times k}\,\zz^{(k)}(t_{n-1})}{k!} }_{\zz(t_{n-1}+\mathbf{c}\,h) + \mathcal{O}(h^{q\I+1})},
\end{split}
\]
which is equivalent to:
\begin{equation}
\label{eqn:Z-external-stages}
\begin{split}
&\left( \mathbf{W}\I \otimes \Idvar \right)\, \Nordsieck_p(h,\zz,t_n) = \left( \mathbf{M}\I(\infty) \,\mathbf{W}\I\otimes \Idvar \right)\, \Nordsieck_p(h,\zz,t_{n-1}) \\
& \qquad + \left( \mathbf{B}\I\,\mathbf{A}\I\,^{-1}\otimes \Idvar \right)\, \zz(t_{n-1}+\mathbf{c}\,h) + \mathcal{O}(h^{q\I+1}).
\end{split}
\end{equation}
This holds in both cases $q\I=p$ and $q\I=p-1$.
Subtracting \eqref{eqn:Z-external-stages} from \eqref{eqn:stiff-stages} and using the stiff stage accuracy equation \eqref{eqn:stiff-stage-global-err} leads to the following recurrence for the stiff global errors \eqref{eqn:global-error-notation}:
\begin{equation}
\label{eqn:stiff-error-recurrence}
\begin{split}
\Delta \zglm^{[n]} &= \Bigl( \mathbf{M}\I(\infty) \otimes \Idvar \Bigr)\,\Delta\zglm^{[n-1]} + \delta^{[n-1]}, \\
 \delta^{[n-1]} &= \left( \mathbf{B}\I\,\mathbf{A}\I\,^{-1} \otimes \Idvar \right)  \,  \Delta Z^{[n]} + \mathcal{O}(h^{q\I+1}) \sim \mathcal{O}(h^{q+1}), 
 \end{split}
\end{equation}
where $q = \min\{ q\E,q\I \}$. The recurrence \eqref{eqn:stiff-error-recurrence} has the solution:
\begin{equation}
\label{eqn:stiff-error-recurrence2}
\Delta \zglm^{[n]} = \bigl(\mathbf{M}\I(\infty)\bigr)^n\,\Delta\zglm^{[0]} + \sum_{k=1}^n  \bigl(\mathbf{M}\I(\infty)\bigr)^{k-1}\,\delta^{[n-k]}.
\end{equation}
Assume that the initialization procedure satisfies
\[
\Delta \zglm^{[0]} \sim \mathcal{O}(h^{p}).
\]
We note that the local truncation errors have the form
\[
\delta^{[k]} = \varphi_q\,h^{q+1}\, \zz^{(q+1)}(t_{k-1}) + \mathcal{O}(h^{q+2}),
\]
and therefore
\[
\delta^{[k]} - \delta^{[k-L]} = \varphi_q\,h^{q+2}\,L\,\zz^{(q+2)}(t_{k-L}) + \mathcal{O}(h^{q+2}) = \mathcal{O}(h^{q+2}).
\]
We now apply Lemma \ref{lem:error-convergence} to obtain the result.
\end{proof}

\begin{corollary}[Order of IMEX-GLM on explicit index-1 DAE]
\label{cor:index-1-order}
Consider an IMEX-GLM method of order $p$ from the class of interest (i.e., with properties enumerated in Definition \ref{def:interesting-methods}). The method applied to the index-1 problem yields a solution that converges with order $p$:
\begin{subequations}
\label{eqn:index-1-full-order}
\begin{align}
\Delta X^{[n]} & \sim \mathcal{O}(h^{p}), &
\Delta\xglm^{[n]} &\sim \mathcal{O}(h^{p}), \\
\Delta Z^{[n]} &\sim \mathcal{O}(h^{p}), &
\Delta\zglm^{[n]} &\sim \mathcal{O}(h^{p}), 
\end{align}
\end{subequations}
\end{corollary}

\begin{remark}[Stiffly accurate GLMs]
If the method is stiffly accurate then $q=p$. Since $\gun(\xglm^{[n]},\zglm^{[n]}) = \gun\big(X^{[n]},Z^{[n]}\big) = \Zero$, the solution obeys the algebraic constraint. 
\end{remark}

\section{Convergence analysis for singular perturbation problems}
\label{sec:SPP-convergence}

\subsection{Solution expansions}
\label{sec:SPP-index-2}

We now formally expand the numerical solutions of method \eqref{eqn:IMEX_SPP} in series of $\varepsilon$:
\begin{equation}
\label{eqn:epsilon-series}
\begin{split}
X^{[n]} = \sum_{k \ge 0} X^{[n],k}\,\varepsilon^k, &\quad
Z^{[n]} = \sum_{k \ge 0} Z^{[n],k}\,\varepsilon^k, \\
\xglm^{[n]} = \sum_{k \ge 0} \xglm^{[n],k}\,\varepsilon^k, &\quad
\zglm^{[n]} = \sum_{k \ge 0} \zglm^{[n],k}\,\varepsilon^k.
\end{split}
\end{equation}
We denote the global errors for each internal and external stage coefficients by:
\begin{equation*}
\begin{split}
\Delta X^{[n],k} &\coloneqq X^{[n],k} - \xx^{k}(t_{n-1} + \mathbf{c} h), \\
\Delta Z^{[n],k} &\coloneqq Z^{[n],k} - \zz^{k}(t_{n-1} + \mathbf{c} h), \\
\Delta\xglm^{[n],k} &\coloneqq \xglm^{[n],k} -  \left(\mathbf{W}\E \otimes \Idvar \right)\,\Nordsieck_p(h,\xx^{k},t_n), \\
\Delta\zglm^{[n],k} &\coloneqq \zglm^{[n],k} -  \left(\mathbf{W}\I \otimes \Idvar \right)\,\Nordsieck_p(h,\zz^{k},t_n).
\end{split}
\end{equation*}

We insert \eqref{eqn:epsilon-series} into \eqref{eqn:IMEX_SPP-intX}, \eqref{eqn:IMEX_SPP-intZ}  and expand the function in Taylor series about the $\mathcal{O}(1)$ terms to obtain:
\begin{subequations}
\label{eqn:IMEX_SPP-epsilon-expansion}
\begin{eqnarray}
\nonumber
\sum_{k \ge 0} X^{[n],k}\,\varepsilon^k &=& h\, \left(\mathbf{A}\E \otimes \Idvar \right) \, \fun(X^{[n],0}\,, Z^{[n],0}) \\
\nonumber
&&+ h\, \left(\mathbf{A}\E \otimes \Idvar \right) \, 
\sum_{\alpha+\beta\ge 1} \frac{\partial^{\alpha+\beta}\fun}{\partial\xx^\alpha\partial\zz^\beta}(X^{[n],0}\,, Z^{[n],0}) \cdot \big( W_{\alpha,\beta}^{[n]} )\quad \\
\label{eqn:IMEX_SPP-intX-epsilon-expansion}
&& + \left(\mathbf{U}\E \otimes \Idvar \right) \, \sum_{k \ge 0} \xglm^{[n-1],k}\,\varepsilon^k, \\
\nonumber
\sum_{k \ge 0} Z^{[n],k}\,\varepsilon^{k+1} &=&  h\, \left(\mathbf{A}\I \otimes \Idvar \right) \, \gun(X^{[n],0}\,, Z^{[n],0}) \\
\nonumber
&&+ h\,\left( \mathbf{A}\I \otimes \Idvar \right) \, \sum_{\alpha+\beta\ge 1} \frac{\partial^{\alpha+\beta}\gun}{\partial\xx^\alpha\partial\zz^\beta}(X^{[n],0}\,, Z^{[n],0}) \cdot \big( W_{\alpha,\beta}^{[n]} )\quad \\
\label{eqn:IMEX_SPP-intZ-epsilon-expansion}
&& + \left(\mathbf{U}\I \otimes \Idvar \right) \,  \sum_{k \ge 0} \zglm^{[n-1],k}\,\varepsilon^{k+1},
\end{eqnarray}
where the function derivatives are applied to the following multinomial arguments:
\begin{equation*}
\begin{split}
\big( W_{\alpha,\beta}^{[n]} \big) &\coloneqq  \Big( \underbrace{\dots, \sum_{m \ge 1} X^{[n],m}\,\varepsilon^m, \dots}_{\alpha\textnormal{ times}}~, ~ \underbrace{\dots,\sum_{\ell \ge 1} Z^{[n],\ell}\,\varepsilon^\ell, \dots}_{\beta\textnormal{ times}} \Big) = \sum_{k \ge 1} \big( W_{\alpha,\beta}^{[n],k} \big)\, \varepsilon^{k}, \\
\big( W_{\alpha,\beta}^{[n],k} \big) &= \sum_{\sum_{i = 1}^\alpha m_i + \sum_{j = 1}^\beta \ell_j = k} 
\Big( X^{[n],m_1},\dots, X^{[n],m_\alpha}, Z^{[n],\ell_1} , \dots, Z^{[n],\ell_\beta} \Big).
\end{split}
\end{equation*}
Similar expressions are obtained by inserting \eqref{eqn:epsilon-series} into the external stage equations \eqref{eqn:IMEX_SPP-extX} and \eqref{eqn:IMEX_SPP-extZ}:
\begin{eqnarray}
\nonumber
\sum_{k \ge 0} \xglm^{[n],k}\,\varepsilon^k &=& h\, \left(\mathbf{B}\E \otimes \Idvar \right) \, \fun(X^{[n],0}\,, Z^{[n],0}) \\ 
\nonumber
&&+ h\, \left(\mathbf{B}\E \otimes \Idvar \right) \, 
\sum_{\alpha+\beta\ge 1} \frac{\partial^{\alpha+\beta}\fun}{\partial\xx^\alpha\partial\zz^\beta}(X^{[n],0}\,, Z^{[n],0}) \cdot \big( W_{\alpha,\beta}^{[n]} )\quad \\
\label{eqn:IMEX_SPP-extX-epsilon-expansion}
&& + \left(\mathbf{V}\E \otimes \Idvar \right) \, \sum_{k \ge 0} \xglm^{[n-1],k}\,\varepsilon^k, \\
\nonumber
\sum_{k \ge 0} \zglm^{[n],k}\,\varepsilon^{k+1} &=& h\, \left(\mathbf{B}\I \otimes \Idvar \right) \, \gun(X^{[n],0}\,, Z^{[n],0}) \\
\nonumber
&&+ h\,\left( \mathbf{B}\I \otimes \Idvar \right) \, \sum_{\alpha+\beta\ge 1} \frac{\partial^{\alpha+\beta}\gun}{\partial\xx^\alpha\partial\zz^\beta}(X^{[n],0}\,, Z^{[n],0}) \cdot \big( W_{\alpha,\beta}^{[n]} )\quad \\
\label{eqn:IMEX_SPP-extZ-epsilon-expansion}
&& + \left(\mathbf{V}\I \otimes \Idvar \right) \,  \sum_{k \ge 0} \zglm^{[n-1],k}\,\varepsilon^{k+1}.
\end{eqnarray}
\end{subequations}

Denote the function derivatives evaluated at the $\mathcal{O}(1)$ terms of the exact solution by:
\begin{equation}
\label{eqn:derivative-notation}
\frac{\partial^{\alpha+\beta}\fun}{\partial\xx^\alpha\partial\zz^\beta}\Big|_{t_{n-1}+ \mathbf{c}h} \coloneqq \frac{\partial^{\alpha+\beta}\fun}{\partial\xx^\alpha\partial\zz^\beta}(\xx^{0}(t_{n-1}+ \mathbf{c}h)\,, \zz^{0}(t_{n-1}+ \mathbf{c}h)),
\end{equation}
and similarly for the derivatives of $\gun$. From \eqref{eqn:index1-stage-global-err} we know that the stage values $(X^{[n],0}\,, Z^{[n],0})$ have global errors $\mathcal{O}(h^{q\E+1})$. Owing to the smoothness of $\fun$, $\gun$, changing the arguments of the function derivatives in \eqref{eqn:IMEX_SPP-epsilon-expansion} from the numerical stage values to the exact solution changes the derivatives by:
\begin{equation}
\label{eqn:derivative-approximation}
\left\Vert \frac{\partial^{\alpha+\beta}\fun}{\partial\xx^\alpha\partial\zz^\beta}(X^{[n],0}\,, Z^{[n],0}) - \frac{\partial^{\alpha+\beta}\fun}{\partial\xx^\alpha\partial\zz^\beta}\Big|_{t_{n-1}+ \mathbf{c}h} \right\Vert \sim \mathcal{O}(h^{q\E+1}),
\end{equation}
and similarly for the derivatives of $\gun$.

Next, we equate powers of $\varepsilon$ in  \eqref{eqn:IMEX_SPP-epsilon-expansion} to obtain equations for each of the expansion coefficients \eqref{eqn:epsilon-series}. 
Using \eqref{eqn:derivative-approximation} to shift the arguments to the exact index-1 DAE solutions the equations for the coefficients of $\varepsilon^k$, $k\ge 1$, are:
\begin{subequations}
\label{eqn:IMEX_SPP-epsilon-k}
\begin{eqnarray}
\nonumber
X^{[n],k} &=&  h\, \left(\mathbf{A}\E \otimes \Idvar \right) \, 
\sum_{1 \le \alpha + \beta\le k} \frac{\partial^{\alpha+\beta}\fun}{\partial\xx^\alpha\partial\zz^\beta}\Big|_{t_{n-1}+ \mathbf{c}h} \cdot \big( W_{\alpha,\beta}^{[n],k} )\quad \\
\label{eqn:IMEX_SPP-intX-epsilon-k}
&& + \left(\mathbf{U}\E \otimes \Idvar \right) \, \xglm^{[n-1],k} + \mathcal{O}(h^{q\E+2}), \\
\nonumber
Z^{[n],k-1} &=& h\,\left( \mathbf{A}\I \otimes \Idvar \right) \,\sum_{1 \le \alpha + \beta\le k} \frac{\partial^{\alpha+\beta}\gun}{\partial\xx^\alpha\partial\zz^\beta}\Big|_{t_{n-1}+ \mathbf{c}h} \cdot \big( W_{\alpha,\beta}^{[n],k} )\quad \\
\label{eqn:IMEX_SPP-intZ-epsilon-k}
&& + \left(\mathbf{U}\I \otimes \Idvar \right) \,  \zglm^{[n-1],k-1}+ \mathcal{O}(h^{q\E+2}),
\end{eqnarray}
\begin{eqnarray}
\nonumber
\xglm^{[n],k} &=& h\, \left(\mathbf{B}\E \otimes \Idvar \right) \, 
\sum_{1 \le \alpha + \beta\le k} \frac{\partial^{\alpha+\beta}\fun}{\partial\xx^\alpha\partial\zz^\beta}\Big|_{t_{n-1}+ \mathbf{c}h} \cdot \big( W_{\alpha,\beta}^{[n],k} )\quad \\
\label{eqn:IMEX_SPP-extX-epsilon-k}
&& + \left(\mathbf{V}\E \otimes \Idvar \right) \, \xglm^{[n-1],k}+ \mathcal{O}(h^{q\E+2}), \\
\nonumber
\zglm^{[n],k-1} &=& h\,\left( \mathbf{B}\I \otimes \Idvar \right) \, \sum_{1 \le \alpha + \beta\le k} \frac{\partial^{\alpha+\beta}\gun}{\partial\xx^\alpha\partial\zz^\beta}\Big|_{t_{n-1}+ \mathbf{c}h} \cdot \big( W_{\alpha,\beta}^{[n],k} )\quad \\
\label{eqn:IMEX_SPP-extZ-epsilon-k}
&& + \left(\mathbf{V}\I \otimes \Idvar \right) \, \zglm^{[n-1],k-1}+ \mathcal{O}(h^{q\E+2}).
\end{eqnarray}
\end{subequations}

Inserting the exact solution $\xx^k$ and its derivative from \eqref{eqn:spp-eps-expansion-k-x} in place of the numerical solutions in \eqref{eqn:IMEX_SPP-intX-epsilon-k} and \eqref{eqn:IMEX_SPP-extX-epsilon-k} leads to residuals $\mathcal{O}(h^{q\E+1})$ and $\mathcal{O}(h^{p+1})$, respectively. Similarly, inserting the exact solution $\zz^{k-1}$ and its derivatives 
\eqref{eqn:spp-eps-expansion-k-z} in place of the numerical solutions in \eqref{eqn:IMEX_SPP-intZ-epsilon-k} and \eqref{eqn:IMEX_SPP-extZ-epsilon-k} leads to residuals $\mathcal{O}(h^{q\I+1})$ and $\mathcal{O}(h^{p+1})$, respectively. Subtracting these relations from \eqref{eqn:IMEX_SPP-epsilon-k} leads to the following recurences of the global errors:
\begin{subequations}
\label{eqn:IMEX_SPP-epsilon-k-globalerr}
\begin{eqnarray}
\nonumber
\Delta X^{[n],k} &=&  h\, \left(\mathbf{A}\E \otimes \Idvar \right) \, 
\Bigl( \fun_\xx\big|_{t_{n-1}+ \mathbf{c}h} \,\Delta X^{[n],k} + \fun_\zz\big|_{t_{n-1}+ \mathbf{c}h} \,\Delta Z^{[n],k} \\
\nonumber
&&\quad + \sum_{2 \le \alpha + \beta\le k} \frac{\partial^{\alpha+\beta}\fun}{\partial\xx^\alpha\partial\zz^\beta}\Big|_{t_{n-1}+ \mathbf{c}h} \cdot \big( W_{\alpha,\beta}^{[n],k} - W_{\alpha,\beta}^{k} \big) \Bigr)\quad \\
\label{eqn:IMEX_SPP-intX-epsilon-k-globalerr}
&& + \left(\mathbf{U}\E \otimes \Idvar \right) \, \Delta\xglm^{[n-1],k} + \mathcal{O}(h^{q+1}), \\
\nonumber
\Delta Z^{[n],k-1} &=& h\,\left( \mathbf{A}\I \otimes \Idvar \right) \,\Bigl( \gun_\xx\big|_{t_{n-1}+ \mathbf{c}h} \,\Delta X^{[n],k} + \gun_\zz\big|_{t_{n-1}+ \mathbf{c}h} \,\Delta Z^{[n],k} \\
\nonumber
&& + \sum_{2 \le \alpha + \beta\le k} \frac{\partial^{\alpha+\beta}\gun}{\partial\xx^\alpha\partial\zz^\beta}\Big|_{t_{n-1}+ \mathbf{c}h} \cdot \big( W_{\alpha,\beta}^{[n],k}-W_{\alpha,\beta}^{k} \big) \Bigr) \quad \\
\label{eqn:IMEX_SPP-intZ-epsilon-k-globalerr}
&& + \left(\mathbf{U}\I \otimes \Idvar \right) \,  \Delta\zglm^{[n-1],k-1}+ \mathcal{O}(h^{q+1}),
\end{eqnarray}
\begin{eqnarray}
\nonumber
\Delta\xglm^{[n],k} &=& h\, \left(\mathbf{B}\E \otimes \Idvar \right) \, 
\Bigl( \fun_\xx\big|_{t_{n-1}+ \mathbf{c}h} \,\Delta X^{[n],k} + \fun_\zz\big|_{t_{n-1}+ \mathbf{c}h} \,\Delta Z^{[n],k} \\
\nonumber
&&+   \sum_{2 \le \alpha + \beta\le k} \frac{\partial^{\alpha+\beta}\fun}{\partial\xx^\alpha\partial\zz^\beta}\Big|_{t_{n-1}+ \mathbf{c}h} \cdot \big( W_{\alpha,\beta}^{[n],k}-W_{\alpha,\beta}^{[k} ) \Bigr)\quad \\
\label{eqn:IMEX_SPP-extX-epsilon-k-globalerr}
&& + \left(\mathbf{V}\E \otimes \Idvar \right) \, \Delta\xglm^{[n-1],k}+ \mathcal{O}(h^{\min(p+1,q+2)}), \\
\nonumber
\Delta\zglm^{[n],k-1} &=& h\,\left( \mathbf{B}\I \otimes \Idvar \right) \, \Bigl( \gun_\xx\big|_{t_{n-1}+ \mathbf{c}h} \,\Delta X^{[n],k} + \gun_\zz\big|_{t_{n-1}+ \mathbf{c}h} \,\Delta Z^{[n],k} \\ 
\nonumber
&&+ \sum_{2 \le \alpha + \beta\le k} \frac{\partial^{\alpha+\beta}\gun}{\partial\xx^\alpha\partial\zz^\beta}\Big|_{t_{n-1}+ \mathbf{c}h} \cdot \big( W_{\alpha,\beta}^{[n],k}-W_{\alpha,\beta}^{k} \big) \Bigr)\quad \\
\label{eqn:IMEX_SPP-extZ-epsilon-k-globalerr}
&& + \left(\mathbf{V}\I \otimes \Idvar \right) \, \Delta\zglm^{[n-1],k-1}+ \mathcal{O}(h^{\min(p+1,q+2)}).
\end{eqnarray}
\end{subequations}
Using the invertibility of $\mathbf{A}\I $ one substitutes \eqref{eqn:IMEX_SPP-intZ-epsilon-k-globalerr} into \eqref{eqn:IMEX_SPP-extZ-epsilon-k-globalerr} to obtain:
\begin{equation}
\label{eqn:IMEX_SPP-substitution}
\begin{split}
& \Delta \zglm^{[n],k-1} = \left( \mathbf{M}\I(\infty) \otimes \Idvar \right) \, \Delta \zglm^{[n-1],k-1} \\
&\qquad + \left( \mathbf{B}\I\mathbf{A}\I\,\!^{-1} \otimes \Idvar \right) \, \Delta Z^{[n],k-1} \\
&\qquad - \left( \mathbf{B}\I \otimes \Idvar \right) \,\sum_{2 \le \alpha + \beta\le k} \frac{\partial^{\alpha+\beta}\gun}{\partial\xx^\alpha\partial\zz^\beta}\Big|_{t_{n-1}+ \mathbf{c}h} \cdot \big( W_{\alpha,\beta}^{[n],k}-W_{\alpha,\beta}^{k} \big)  \\
&\qquad + \mathcal{O}(h^{q+1}).
\end{split}
\end{equation}
Note that the entries of the multilinear argument difference $W_{\alpha,\beta}^{[n],k}-W_{\alpha,\beta}^{k}$ for $\alpha+\beta \le k$ contain global errors of stage values $\Delta X^{[n],j}$, $\Delta Z^{[n],j}$ up to index $j \le k-1$. Consequently:
\begin{equation}
\label{eqn:multilinear-bound}
\big\Vert W_{\alpha,\beta}^{[n],k}-W_{\alpha,\beta}^{k} \big\Vert = \max_{0 \le j \le k-1} \big(\Vert \Delta X^{[n],j}\Vert ,\Vert Z^{[n],j}\Vert \big).
\end{equation}

\begin{remark}[GLM solutions of DAEs]
\label{rem:GLM-applied_to-index_k}
Let
\begin{equation}
\label{eqn:GLM_k_ell}
\begin{split}
\kappa^{[n]} &\coloneqq \fun\left( X^{[n]}, Z^{[n]} \right) = \fun\left( \sum_{k \ge 0} X^{[n],k}\,\varepsilon^k, \sum_{k \ge 0} Z^{[n],k}\,\varepsilon^k \right)
= \sum_{k \ge 0} \kappa^{[n],k}\,\varepsilon^k, \\
\kappa^{[n],k} &= 
\sum_{1 \le \alpha + \beta\le k} \frac{\partial^{\alpha+\beta}\fun}{\partial\xx^\alpha\partial\zz^\beta}\left( X^{[n],0}, Z^{[n],0} \right) \cdot \big( W_{\alpha,\beta}^{[n],k} ), \\
\varepsilon\,\ell^{[n]} &\coloneqq \gun\left( X^{[n]}, Z^{[n]} \right) = \gun\left( \sum_{k \ge 0} X^{[n],k}\,\varepsilon^k, \sum_{k \ge 0} Z^{[n],k}\,\varepsilon^k \right)
= \sum_{k \ge 0} \ell^{[n],k}\,\varepsilon^{k+1}, \\
\ell^{[n],k-1} &=  
\sum_{1 \le \alpha + \beta\le k} \frac{\partial^{\alpha+\beta}\gun}{\partial\xx^\alpha\partial\zz^\beta}\left( X^{[n],0}, Z^{[n],0} \right) \cdot \big( W_{\alpha,\beta}^{[n],k} )
\end{split}
\end{equation}
We see that \eqref{eqn:IMEX_SPP-intX-epsilon-k} is the result of applying the IMEX GLM to solve the index-$(k+1)$ differential algebraic equation \eqref{eqn:spp-eps-expansion-k-x} for the unknown functions $\xx^k$ and $\zz^{k-1}$.
\begin{equation}
\label{eqn:GLM-applied_to-index_k}
\begin{split}
X^{[n],k} &=  h\, \left(\mathbf{A}\E \otimes \Idvar \right) \,\kappa^{[n],k}  + \left(\mathbf{U}\E \otimes \Idvar \right) \, \xglm^{[n-1],k}, \\
Z^{[n],k-1} &= h\,\left( \mathbf{A}\I \otimes \Idvar \right) \,\ell^{[n],k-1} + \left(\mathbf{U}\I \otimes \Idvar \right) \,  \zglm^{[n-1],k-1}, \\
\xglm^{[n],k} &= h\, \left(\mathbf{B}\E \otimes \Idvar \right) \, \kappa^{[n],k} + \left(\mathbf{V}\E \otimes \Idvar \right) \, \xglm^{[n-1],k}, \\
\zglm^{[n],k-1} &= h\,\left( \mathbf{B}\I \otimes \Idvar \right) \,\ell^{[n],k-1}  + \left(\mathbf{V}\I \otimes \Idvar \right) \, \zglm^{[n-1],k-1}.
\end{split}
\end{equation}
\end{remark}

\subsection{The $\varepsilon^0$ coefficients}

From \eqref{eqn:IMEX_SPP-epsilon-expansion} we see that the terms $X^{[n],0}$, $Z^{[n],0}$, $\xglm^{[n],0}$, and $\zglm^{[n],0}$ are the numerical solutions of the index-1 DAE system, which have been already discussed in Theorem \ref{thm:index-k-order},
The global errors $\Delta X^{[n],0}$, $\Delta Z^{[n],0}$, $\Delta \xglm^{[n],0}$, and $\Delta \zglm^{[n],0}$ are given by equations \eqref{eqn:index-1-order}.

\subsection{The $\varepsilon^1$ coefficients}

We next equate the terms in $\varepsilon^1$ in \eqref{eqn:IMEX_SPP-epsilon-expansion}. Specifically, the global error recurrences \eqref{eqn:IMEX_SPP-epsilon-k-globalerr} for $k=1$ are:
\begin{subequations}
\label{eqn:IMEX_SPP-index2-err}
\begin{eqnarray}
\label{eqn:IMEX_SPP-intX-index2-err}
\Delta X^{[n],1} &=& h\, \left(\mathbf{A}\E \otimes \Idvar \right) \, \left( \fun_\xx\big|_{t_{n-1}+\mathbf{c}h}\,\Delta X^{[n],1} + \fun_\zz\big|_{t_{n-1}+\mathbf{c}h}\,\Delta Z^{[n],1} \right) \\
\nonumber
&& + \left(\mathbf{U}\E \otimes \Idvar \right) \, \Delta \xglm^{[n-1],1} + \mathcal{O}(h^{q+1}),\qquad \\ 
\label{eqn:IMEX_SPP-intZ-index2-err}
\underbrace{\Delta Z^{[n],0}}_{\mathcal{O}(h^{q+1})} &=& h\,\left( \mathbf{A}\I \otimes \Idvar \right) \, \left( \gun_\xx\big|_{t_{n-1}+\mathbf{c}h}\,\Delta X^{[n],1} + \gun_\zz\big|_{t_{n-1}+\mathbf{c}h}\,\Delta Z^{[n],1} \right)  \\
\nonumber
&& + \left(\mathbf{U}\I \otimes \Idvar \right) \,\underbrace{\Delta \zglm^{[n-1],0}}_{\mathcal{O}(h^{\nu})} + \mathcal{O}(h^{q+1}),\qquad \\ 
\label{eqn:IMEX_SPP-extX-index2-err}
\Delta \xglm^{[n],1} &=& h\, \left(\mathbf{B}\E \otimes \Idvar \right) \, \left( \fun_\xx\big|_{t_{n-1}+\mathbf{c}h}\,\Delta X^{[n],1} + \fun_\zz\big|_{t_{n-1}+\mathbf{c}h}\,\Delta Z^{[n],1} \right)   \\
\nonumber
&& + \left(\mathbf{V}\E \otimes \Idvar \right) \, \Delta \xglm^{[n-1],1} + \mathcal{O}(h^{\min(p+1,q+2)}), \\
\label{eqn:IMEX_SPP-extZ-index2-err}
 \underbrace{\Delta \zglm^{[n],0}}_{\mathcal{O}(h^{\nu})} &=& h\,\left( \mathbf{B}\I \otimes \Idvar \right)  \,  \left( \gun_\xx\big|_{t_{n-1}+\mathbf{c}h}\,\Delta X^{[n],1} + \gun_\zz\big|_{t_{n-1}+\mathbf{c}h}\,\Delta Z^{[n],1} \right)  \\
\nonumber
&& + \left(\mathbf{V}\I \otimes \Idvar \right) \, \underbrace{\Delta \zglm^{[n-1],0}}_{\mathcal{O}(h^{\nu})} + \mathcal{O}(h^{\min(p+1,q+2)}).\qquad
\end{eqnarray}
\end{subequations}
From \eqref{eqn:IMEX_SPP-intZ-index2-err} we have that:
\begin{equation}
\label{eqn:dz1-error}
\Delta Z^{[n],1} = - (\gun_\zz^{-1} \gun_\xx)\big|_{t_{n-1}+\mathbf{c}h}\,\Delta X^{[n],1}  + \mathcal{O}(h^{\min(\nu-1,q)}),
\end{equation}
and inserting this solution into \eqref{eqn:IMEX_SPP-intX-index2-err} and \eqref{eqn:IMEX_SPP-extX-index2-err} gives the following the global error relations:
\begin{eqnarray*}
\Delta X^{[n],1} &=& h\, \left(\mathbf{A}\E \otimes \Idvar \right) \, \left( (\fun_\xx - \fun_\zz\,\gun_\zz^{-1}\,\gun_\xx)\big|_{t_{n-1}+\mathbf{c}h}\,\Delta X^{[n],1}  \right) \\
\nonumber
&& + \left(\mathbf{U}\E \otimes \Idvar \right) \, \Delta \xglm^{[n-1],1} + \mathcal{O}(h^{\min(\nu,q+1)}),\qquad \\ 
\Delta \xglm^{[n],1} &=& h\, \left(\mathbf{B}\E \otimes \Idvar \right) \, \left( (\fun_\xx - \fun_\zz\,\gun_\zz^{-1}\,\gun_\xx)\big|_{t_{n-1}+\mathbf{c}h}\,\Delta X^{[n],1}  \right)   \\
&& + \left(\mathbf{V}\E \otimes \Idvar \right) \, \Delta \xglm^{[n-1],1} + \mathcal{O}(h^{\min(\nu,q+1)}) \\
&=& \mathbf{M}\E \left( h\,(\fun_\xx - \fun_\zz\,\gun_\zz^{-1}\,\gun_\xx)\big|_{t_{n-1}+\mathbf{c}h} \right) \,  \Delta \xglm^{[n-1],1} + \mathcal{O}(h^{\min(\nu,q+1)}),
\end{eqnarray*}
where $\mathbf{M}\E$ is the stability matrix of the explicit GLM component. Using the stability of $\mathbf{M}\E$ and the convergence of the explicit scheme we have:
\begin{subequations}
\label{eqn:index2-order}
\begin{equation}
\label{eqn:x-index2-order}
\begin{split}
\Delta \xglm^{[n],1} &\sim \mathcal{O}\bigl( h^{\min(\nu-1,q)} \bigr), \\
\Delta X^{[n],1} &\sim \mathcal{O}\bigl( h^{\min(\nu-1,q)} \bigr).
\end{split}
\end{equation}
Consequently, \eqref{eqn:dz1-error} and the invertibility of $\gun_\zz$ reveal that
\begin{equation}
\label{eqn:z-index2-order}
\begin{split}
\Delta Z^{[n],1} &\sim \mathcal{O}\bigl( h^{\min(\nu-1,q)} \bigr).
\end{split}
\end{equation}
From \eqref{eqn:IMEX_SPP-substitution} with $k=2$, using the error bound \eqref{eqn:multilinear-bound}, we have:
\begin{equation}
\Delta\zglm^{[n],1} = \left(\mathbf{M}\I(\infty) \otimes \Idvar \right)  \, \Delta\zglm^{[n-1],1} + \mathcal{O}\bigl( h^{\min(\nu-1,q)} \bigr),
\end{equation}
\end{subequations}
and solving this recurrence gives the order of $\Delta\zglm^{[n],1}$.

\begin{theorem}[Order of IMEX-GLM on differential-algebraic problems]
\label{thm:index-k-order}
Consider the singular perturbation problem \eqref{eqn:SPP}, with an initial solution \eqref{eqn:SPP-epsilon-expansion} that is consistent with respect to \eqref{eqn:dae-eps-0}, \eqref{eqn:expansion-eps-nu}. Consider an IMEX-GLM method of order $p$ from the class of interest (i.e., with properties enumerated in Definition \ref{def:interesting-methods}). 

Application of this IMEX GLM to solve the index-$(k+1)$ differential algebraic equations \eqref{eqn:spp-eps-expansion-k-x} 
yields solutions with the following global errors: 
\begin{align*}
\Delta X^{[n]} &= \sum_{k \ge 0} \Delta X^{[n],k}\, \varepsilon^k, &
\Delta\xglm^{[n]} &= \sum_{k \ge 0} \Delta\xglm^{[n],k}\, \varepsilon^k, \\
\Delta Z^{[n]} &= \sum_{k \ge 0} \Delta Z^{[n],k}\, \varepsilon^k, &
\Delta\zglm^{[n]} &= \sum_{k \ge 0} \Delta\zglm^{[n],k}\, \varepsilon^k, 
\end{align*}
where, for $0 \le k \le p$:
\begin{subequations}
\label{eqn:index-k-global-errors}
\begin{align}
\Delta X^{[n],k} &\sim \mathcal{O}\big( h^{p-k} \big), &
\Delta\xglm^{[n],k} & \sim \mathcal{O}\big( h^{p-k} \big), \\
\Delta Z^{[n],k} & \sim \mathcal{O}\big( h^{p-k} \big), &
\Delta\zglm^{[n],k} & \sim \mathcal{O}\big( h^{p-k} \big). 
\end{align}
\end{subequations}
\end{theorem}

\begin{proof}
Equations \eqref{eqn:index-1-full-order} prove the case $k=0$. Equations \eqref{eqn:index2-order} prove the case $k=1$. The rest follows by induction on $k$.
\end{proof}

 
\begin{remark}[Traditional implicit GLMs]
Theorem 
Hairer and collaborators \cite{Hairer_1988_RK-SPP}, \cite[Chapter VI.3]{Hairer_book_II} (for Runge-Kutta methods) and Schneider  \cite{Schneider_1993_GLM-SPP}  (for GLMs) consider the quantities:
\begin{equation*}
\begin{split}
& \Delta X^{[n],1} - (\fun_z\,\gun_z^{-1})\big|_{t_{n-1+\mathbf{c} h}}\,\Delta Z^{[n],0} , \\
& \mathbf{u}^{[n],1} \coloneqq \Delta \xglm^{[n],1} - (\fun_z\,\gun_z^{-1})\big|_{t_n}\,\Delta \zglm^{[n],0} = \Delta \xglm^{[n],1} + \mathcal{O}(h^\nu).
\end{split}
\end{equation*}
and study the convergence of the latter. 
In case of Runge-Kutta methods, and in case of traditional implicit GLMs where $\mathbf{A}\E \leftarrow \mathbf{A}\I$, $\mathbf{B}\E \leftarrow \mathbf{B}\I$,  $\mathbf{U}\E \leftarrow \mathbf{U}\I$, and  $\mathbf{V}\E \leftarrow \mathbf{V}\I$, one can write a recurrence for $\mathbf{u}^{[n],1}$ by taking linear combinations of the internal and external stage equations for $\xx$ and $\zz$ variables. Solving this recurrence shows that $\mathbf{u}^{[n],1} \sim \mathcal{O}(h^\nu)$, and therefore $\Delta \xglm^{[n],1} \sim \mathcal{O}(h^\nu)$. By induction, the asymptotic orders of $\Delta \xglm^{[n],k}$ are one order higher than those of $\Delta \zglm^{[n],k}$.

For general IMEX GLMs where the explicit and the implicit method coefficients are different such a recurrence cannot be formed. This is due to the fact that each of the terms $\Delta \xglm^{[n],1}$ and $(\fun_z\,\gun_z^{-1})\big|_{t_n}\,\Delta \zglm^{[n],0}$ needs to carry different scaling coefficients when forming the linear combinations of internal stages and of external stages. 
\end{remark}

\begin{theorem}[Order of IMEX-GLM on singular perturbation problem]
\label{thm:spp-order}
Consider the singular perturbation problem \eqref{eqn:SPP}, with an initial solution \eqref{eqn:SPP-epsilon-expansion} that is consistent with respect to \eqref{eqn:dae-eps-0}, \eqref{eqn:expansion-eps-nu}. Consider an IMEX-GLM method of order $p$ from the class of interest (i.e., with properties enumerated in Definition \ref{def:interesting-methods}). 
In addition we assume that the stability matrix of the implicit component has a spectral radius bounded away from one in an entire slab of the complex plane:
\begin{equation}
\label{eqn:eigenvalue-assumption}
\rho\left( \mathbf{M}\I(z) \right) \le \alpha < 1 \quad \forall z \in \Co \,:\, Re\,z \le -D.
\end{equation}
Application of this IMEX GLM to solve the singular perturbation problem \eqref{eqn:SPP}  yields solutions with the following global errors: 
\begin{subequations}
\label{eqn:spp-order}
\begin{align}
\Delta\xglm^{[n]} &= \sum_{k = 0}^p\, \underbrace{\Delta\xglm^{[n],k}}_{\mathcal{O}\big( h^{p-k} \big)}\, \varepsilon^k + \Delta_\xglm^{[n]}, \qquad  \Delta_\xglm^{[n]} \sim \mathcal{O}(\varepsilon^{p+1}/h), \\
\Delta\zglm^{[n]} &= \sum_{k = 0}^p\, \underbrace{\Delta\zglm^{[n],k}}_{\mathcal{O}\big( h^{p-k} \big)}\, \varepsilon^k + \Delta_\zglm^{[n]}, \qquad  \Delta_\zglm^{[n]} \sim  \mathcal{O}(\varepsilon^{p+1}/h), 
\end{align}
\end{subequations}
where the size of the error coefficients $\Delta\xglm^{[n],k}$, $\Delta\zglm^{[n],k}$ for $0 \le k \le p$ is given by \eqref{eqn:index-k-global-errors}. The estimates hold uniformly for $h \le h_0$, $n h \le Const$, and $\varepsilon \le D\,h$.
\end{theorem}

\begin{proof}

The proof of Theorem \ref{thm:spp-order} is divided in several steps. The reasoning parallels very closely the proof for Runge-Kutta methods \cite{Hairer_1988_RK-SPP}, \cite[Theorem VI.3.8]{Hairer_book_II}, and the one for classical GLMs \cite{Schneider_1993_GLM-SPP}.

\paragraph{Technical approach.}

For each of the numerical solution variables define the truncated $\varepsilon$-series and the differences between the numerical solutions and these series:
\begin{equation}
\label{eqn:truncated_series}
\begin{split}
\widehat{X}^{[n]} \coloneqq \sum_{k = 0}^{p} X^{[n],k}\,\varepsilon^k, \quad  
\bigtriangleup  X^{[n]} \coloneqq X^{[n]} - \widehat{X}^{[n]},
\\
\widehat{\xglm}^{[n]} \coloneqq \sum_{k = 0}^{p} \xglm^{[n],k}\,\varepsilon^k, \quad  
\bigtriangleup \xglm^{[n]} \coloneqq \xglm^{[n]}  - \widehat{\xglm}^{[n]},
\\
\widehat{\kappa}^{[n]} \coloneqq \sum_{k = 0}^{p} \kappa^{[n],k}\,\varepsilon^k, \quad  
\bigtriangleup \kappa^{[n]} \coloneqq \kappa^{[n]} - \widehat{\kappa}^{[n]}, 
\end{split}
\end{equation}
and similarly for $\widehat{Z}^{[n]}$, $\widehat{\zglm}^{[n]}$, $\widehat{\ell}^{[n]}$
The global errors have the form:
\begin{subequations}
\label{eqn:truncated_global_errors}
\begin{equation}
\begin{split}
\xglm^{[n]} - \mathbf{W}\E\, \eta_p(t_n,h,\xx) &= \sum_{k = 0}^{p} \Delta \xglm^{[n],k}\,\varepsilon^k + \Delta_\xglm^{[n]}, \\
\Delta_\xglm^{[n]} &= \underbrace{\sum_{k \ge p+1} \xglm^{[n],k} \,\varepsilon^k}_{\xglm^{[n]}-\widehat{\xglm}^{[n]}} - 
\underbrace{\sum_{k \ge p+1} \mathbf{W}\E\, \eta_p(t_n,h,\xx^k) \,\varepsilon^k}_{\mathcal{O}(\varepsilon^{p+1})}, 
\end{split}
\end{equation}
where for the second term we used the known fact that the coefficients of $\varepsilon^k$ in the $\varepsilon$-expansion of the exact solution $\xx$ are smooth, and therefore uniformly bounded, for any $k$. Consequently:
\begin{equation}
\label{eqn:remainder-estimate}
\begin{split}
\Vert \Delta_\xglm^{[n]} \Vert &= \Vert \bigtriangleup \xglm^{[n]} \Vert + \mathcal{O}(\varepsilon^{p+1}), \\
\Vert \Delta_\zglm^{[n]} \Vert &= \Vert \bigtriangleup  \zglm^{[n]} \Vert + \mathcal{O}(\varepsilon^{p+1}).
\end{split}
\end{equation}
\end{subequations}
In order to obtain bounds on the remainder terms $\Delta_\xglm^{[n]}$, $\Delta_\zglm^{[n]}$ we look to bund the solution differences between the numerical solution and its truncated version $\bigtriangleup  \xglm^{[n]}$ and $\bigtriangleup  \zglm^{[n]}$.

\paragraph{Step 1. Estimate differences between the internal stages.}
From \eqref{eqn:GLM_k_ell} we have:
\begin{subequations}
\label{eqn:hat-approximations}
\begin{equation}
\label{eqn:hat-approximation-f}
\begin{split}
\fun\left( \widehat{X}^{[n]}, \widehat{Z}^{[n]} \right) &= \sum_{k=0}^p \varepsilon^k\,
\sum_{1 \le \alpha + \beta\le k} \frac{\partial^{\alpha+\beta}\fun}{\partial\xx^\alpha\partial\zz^\beta}\left( X^{[n],0}, Z^{[n],0} \right) \cdot \big( W_{\alpha,\beta}^{[n],k} ) \\
&+ \varepsilon^{p+1}\,\sum_{1 \le \alpha + \beta\le p} \frac{\partial^{\alpha+\beta}\fun}{\partial\xx^\alpha\partial\zz^\beta}\left( X^{[n],0}, Z^{[n],0} \right) \cdot \big( W_{\alpha,\beta}^{[n],p+1} ) \\
&= \widehat{\kappa}^{[n]} + \mathcal{O}(\varepsilon^{p+1}),
\end{split}
\end{equation}
where we use the fact that, and that $W_{\alpha,\beta}^{[n],k}$ uses only arguments  $\widehat{X}^{[n],j}=X^{[n],j}$ with $j < k$. Similarly:
\begin{equation}
\label{eqn:hat-approximation-g}
\begin{split}
\gun\left( \widehat{X}^{[n]}, \widehat{Z}^{[n]} \right) &= \sum_{k=0}^p \varepsilon^k\,
\sum_{1 \le \alpha + \beta\le k} \frac{\partial^{\alpha+\beta}\gun}{\partial\xx^\alpha\partial\zz^\beta}\left( X^{[n],0}, Z^{[n],0} \right) \cdot \big( W_{\alpha,\beta}^{[n],k} ) \\
&+ \varepsilon^{p+1}\,\sum_{1 \le \alpha + \beta\le p} \frac{\partial^{\alpha+\beta}\gun}{\partial\xx^\alpha\partial\zz^\beta}\left( X^{[n],0}, Z^{[n],0} \right) \cdot \big( W_{\alpha,\beta}^{[n],p+1} ) \\
&= \varepsilon\,\widehat{\ell}^{[n]} + \varepsilon^{p+1} \ell^{[n],p} \\
&= \varepsilon\,\widehat{\ell}^{[n]} + \mathcal{O}(\varepsilon^{p+1}/h).
\end{split}
\end{equation}
The last equality come from \eqref{eqn:GLM-applied_to-index_k} and \eqref{eqn:index-k-global-errors} (considering the $\mathcal{O}(1)$ global error with $k=p$) which give:
\begin{equation}
\begin{split}
h\,\ell^{[n],p} &=  \left( \mathbf{A}\I\,\!^{-1} \otimes \Idvar \right) \,\left(Z^{[n],p} - \left(\mathbf{U}\I \otimes \Idvar \right) \,  \zglm^{[n-1],p}\right) \\
&=  \left( \mathbf{A}\I\,\!^{-1} \otimes \Idvar \right) \,\left(\zz^{p}(t_{n-1}+\mathbf{c} h) - \zz^p(t_{n-1}) + \mathcal{O}(h) + \mathcal{O}(1) \right) \\
& \Rightarrow \quad \ell^{[n],p} \sim \mathcal{O}(1/h).
\end{split}
\end{equation}
\end{subequations}

By adding the solution equations \eqref{eqn:GLM-applied_to-index_k} for $k=0,\dots,p$, after scaling each component by the corresponding $\varepsilon^k$, leads to:
\begin{equation}
\label{eqn:hat-reccurence}
\begin{split}
&\widehat{X}^{[n]} - \left(\mathbf{U}\E \otimes \Idvar \right) \, \widehat{\xglm}^{[n-1]} =  h\, \left(\mathbf{A}\E \otimes \Idvar \right) \,\widehat{\kappa}^{[n]} \\
&\qquad = h\,\left( \mathbf{A}\E \otimes \Idvar \right) \,\fun\left( \widehat{X}^{[n]}, \widehat{Z}^{[n]}\right) + \mathcal{O}(h\varepsilon^{p+1}), \\
& \varepsilon\,\left(\widehat{Z}^{[n]} - \left(\mathbf{U}\I \otimes \Idvar \right) \, \widehat{\zglm}^{[n-1]} \right) = h\,\left( \mathbf{A}\I \otimes \Idvar \right) \,\varepsilon\,\widehat{\ell}^{[n]}  \\
&\qquad = h\,\left( \mathbf{A}\I \otimes \Idvar \right) \,\gun\left( \widehat{X}^{[n]}, \widehat{Z}^{[n]}\right) + \mathcal{O}(\varepsilon^{p+1}).
\end{split}
\end{equation}
where the last equalities are obtained form the approximations \eqref{eqn:hat-approximations}. 

Assume the initial condition satisfy \eqref{eqn:perturbed-initial-condition} (a property that will be justified later by induction). Application of Lemma \ref{lem:perturbations-on-internal-stages} to the perturbed scheme \eqref{eqn:hat-reccurence} gives the following relation for internal stage differences:
\begin{equation*}
\begin{split}
\bigtriangleup X^{[n]} &=   C_1^{[n]}\,\bigtriangleup \xglm^{[n-1]} + \varepsilon\,C_2^{[n]}\, \bigtriangleup \zglm^{[n-1]}  + \mathcal{O}(\varepsilon^{p+1}/h), \\ 
\bigtriangleup Z^{[n]} &=   C_3^{[n]}\,\bigtriangleup \xglm^{[n-1]} + \frac{\varepsilon}{h}\, C_4^{[n]}\,\bigtriangleup \zglm^{[n-1]}  + \mathcal{O}(\varepsilon^{p+1}/h).
\end{split}
\end{equation*}
Since $C_i^{[n]}$ are uniformly bounded matrices we will use (with a slight abuse) the following notation:
\begin{equation}
\label{eqn:difference-internal-stages}
\begin{split}
\bigtriangleup X^{[n]} &=   \mathcal{O}(1)\,\bigtriangleup \xglm^{[n-1]} + \mathcal{O}(\varepsilon)\, \bigtriangleup \zglm^{[n-1]}  + \mathcal{O}(\varepsilon^{p+1}/h), \\ 
\bigtriangleup Z^{[n]} &=   \mathcal{O}(1)\,\bigtriangleup \xglm^{[n-1]} + \mathcal{O}\left(\varepsilon/h\right)\,\bigtriangleup \zglm^{[n-1]}  + \mathcal{O}(\varepsilon^{p+1}/h).
\end{split}
\end{equation}


\paragraph{Step 2. Estimate differences between the external stages.}
From \eqref{eqn:GLM-applied_to-index_k} and \eqref{eqn:hat-approximation-f} it follows that:
\begin{equation*}
\begin{split}
\widehat{\xglm}^{[n]} &= h\, \left(\mathbf{B}\E \otimes \Idvar \right) \, \widehat{\kappa}^{[n]} + \left(\mathbf{V}\E \otimes \Idvar \right) \, \widehat{\xglm}^{[n-1]} \\
&= h\, \left(\mathbf{B}\E \otimes \Idvar \right) \, \fun\left( \widehat{X}^{[n]}, \widehat{Z}^{[n]} \right) + \left(\mathbf{V}\E \otimes \Idvar \right) \, \widehat{\xglm}^{[n-1]} + \mathcal{O}(h\,\varepsilon^{p+1}), \\
\bigtriangleup\xglm^{[n]}  &= h\, \left(\mathbf{B}\E \otimes \Idvar \right) \, \left( \fun\left( {X}^{[n]}, {Z}^{[n]} \right) - \fun\left( \widehat{X}^{[n]}, \widehat{Z}^{[n]} \right) \right) \\
&\quad + \left(\mathbf{V}\E \otimes \Idvar \right) \,  \bigtriangleup\xglm^{[n-1]} + \mathcal{O}(h\,\varepsilon^{p+1}) \\
&= h\, \left(\mathbf{B}\E \otimes \Idvar \right) \, \left( F_X^{[n]}\, \bigtriangleup X^{[n]} + F_Z^{[n]}\,  \bigtriangleup{Z}^{[n]} \right) \\
&\quad + \left(\mathbf{V}\E \otimes \Idvar \right) \,  \bigtriangleup\xglm^{[n-1]} + \mathcal{O}(h\,\varepsilon^{p+1}),
\end{split}
\end{equation*}
where $F_X^{[n]}$, $F_Z^{[n]}$ are Jacobians resulting from the mean value theorem that are uniformly bounded over the time interval of interest.
%
%
Using \eqref{eqn:difference-internal-stages} yields:
\begin{equation}
\label{eqn:dx-bound}
\begin{split}
\bigtriangleup\xglm^{[n]} 
&=  \left(\mathbf{V}\E \otimes \Idvar + \mathcal{O}(h) \right)\,\bigtriangleup\xglm^{[n-1]} + \mathcal{O}\left(\varepsilon\right)\,\bigtriangleup \zglm^{[n-1]}  + \mathcal{O}(\varepsilon^{p+1}/h).
\end{split}
\end{equation}

From \eqref{eqn:GLM-applied_to-index_k} and the difference bound \eqref{eqn:difference-internal-stages} it results that:
\begin{equation*}
\begin{split}
& \varepsilon\,\bigtriangleup\ell^{[n]} =  \gun\left( {X}^{[n]}, {Z}^{[n]} \right) - \gun\left( \widehat{X}^{[n]}, \widehat{Z}^{[n]} \right) + \mathcal{O}(\varepsilon^{p+1}/h) \\
&\quad = G_X^{[n]}\,\bigtriangleup{X}^{[n]} + G_Z^{[n]}\,\bigtriangleup{Z}^{[n]} +  \mathcal{O}(\varepsilon^{p+1}/h) \\
&\quad = (\Idstage \otimes \gun_\zz\big|_{t_n})\,\bigtriangleup{Z}^{[n]} + \mathcal{O}\left( h\right)\, \bigtriangleup Z^{[n]}  
+ \mathcal{O}\left( 1 \right)\, \bigtriangleup X^{[n]}  + \mathcal{O}(\varepsilon^{p+1}/h),
\end{split}
\end{equation*}
and using \eqref{eqn:difference-internal-stages}:
\begin{equation*}
 \varepsilon\,\bigtriangleup\ell^{[n]} = (\Idstage \otimes \gun_\zz\big|_{t_n})\,\bigtriangleup{Z}^{[n]} +
  \mathcal{O}(1)\,\bigtriangleup \xglm^{[n-1]} + \mathcal{O}(\varepsilon)\, \bigtriangleup \zglm^{[n-1]}  + \mathcal{O}(\varepsilon^{p+1}/h).
\end{equation*}
Inserting
\begin{equation*}
\begin{split}
\bigtriangleup{Z}^{[n]}  &=( \mathbf{U}\I \otimes \Idvar) \, \bigtriangleup\zglm^{[n-1]}  + h\,(\mathbf{A}\I \otimes \Idvar)\, \bigtriangleup\ell^{[n]}
\end{split}
\end{equation*}
leads to:
\begin{equation*}
\begin{split}
& \varepsilon\,\bigtriangleup\ell^{[n]}  = 
h\,(\mathbf{A}\I \otimes \gun_\zz\big|_{t_n})\,\bigtriangleup\ell^{[n]} 
+ ( \mathbf{U}\I \otimes \gun_\zz\big|_{t_n}) \, \bigtriangleup\zglm^{[n-1]}   \\
&\quad +  \mathcal{O}(1)\,\bigtriangleup \xglm^{[n-1]} + \mathcal{O}(\varepsilon)\, \bigtriangleup \zglm^{[n-1]}  + \mathcal{O}(\varepsilon^{p+1}/h) \\
& h\,\bigtriangleup\ell^{[n]} = \left( (\varepsilon/h) \Id - \,(\mathbf{A}\I \otimes \gun_\zz\big|_{t_n}) \right)^{-1} ( \mathbf{U}\I \otimes \gun_\zz\big|_{t_n}) \,  \bigtriangleup\zglm^{[n-1]}  \\
&\quad +  \mathcal{O}(1)\,\bigtriangleup \xglm^{[n-1]} + \mathcal{O}(\varepsilon)\, \bigtriangleup \zglm^{[n-1]}  + \mathcal{O}(\varepsilon^{p+1}/h).
\end{split}
\end{equation*}
where we used the fact that $(\varepsilon/h)\, \Id_{\nvar s} - (\mathbf{A}\I \otimes \gun_\zz\big|_{t_n}) $ has a bounded inverse (apply Theorem \ref{thm:vonNeumann-matrix} with $\mathbf{G} = \Id_{\nvar s}$, and use the invertibility of $\mathbf{A}\I$ and the negative logarithmic norm $\mu(\gun_\zz) \le -1$).

Use now the external stage equation we obtain:
\begin{equation}
\label{eqn:dz-bound}
\begin{split}
\bigtriangleup\zglm^{[n]} &= \left( \mathbf{B}\I \otimes \Idvar \right) \,h\,\bigtriangleup\ell^{[n]}   + \left(\mathbf{V}\I \otimes \Idvar \right) \, \bigtriangleup\zglm^{[n-1]} \\
&= \mathbf{M}\I\left( (h/\varepsilon)\,\gun_\zz\big|_{t_n} \right)\,\bigtriangleup\zglm^{[n-1]}\\
&\quad +  \mathcal{O}(1)\,\bigtriangleup \xglm^{[n-1]} + \mathcal{O}(\varepsilon)\, \bigtriangleup \zglm^{[n-1]}  + \mathcal{O}(\varepsilon^{p+1}/h).
\end{split}
\end{equation}

\paragraph{Step 3. Solve recurrence for the external stages.}

Equations \eqref{eqn:dx-bound} and \eqref{eqn:dz-bound} lead to the recurrence:
\begin{equation}
\label{eqn:bound}
\begin{split}
\begin{bmatrix} \bigtriangleup\xglm^{[n]} \\ \bigtriangleup\zglm^{[n]} \end{bmatrix}
&=  \begin{bmatrix} \mathbf{V}\E \otimes \Idvar + \mathcal{O}(h) & \mathcal{O}\left(\varepsilon\right) \\ \mathcal{O}(1)  & \mathbf{M}\I\left( (h/\varepsilon)\,\gun_\zz\big|_{t_n} \right) + \mathcal{O}(\varepsilon) \end{bmatrix}
\begin{bmatrix} \bigtriangleup\xglm^{[n-1]} \\ \bigtriangleup\zglm^{[n-1]} \end{bmatrix} \\
&\qquad + 
\begin{bmatrix} \mathcal{O}(\varepsilon^{p+1}) \\ \mathcal{O}(\varepsilon^{p+1}/h) \end{bmatrix}.
\end{split}
\end{equation}
%
%
Since $\mu(\gun_\zz\big|_{t}) \le -1$ it follows from Lemma \ref{lem:eig-log-norm} that the eigenvalues $\lambda_i$ of $(h/\varepsilon)\,\gun_\zz\big|_{t}$ have real parts smaller than or equal to $-h/\varepsilon$. Therefore, using assumption \eqref{eqn:eigenvalue-assumption},
\[
\rho\left( \mathbf{M}\I\left( (h/\varepsilon)\,\gun_\zz\big|_{t} \right) \right) \le \max_{i}\,\rho\left( \mathbf{M}\I\left( \lambda_i \right) \right) \le \alpha < 1.
\]
Apply Lemma \ref{lem:perturbations-on-external-stages} to obtain the bounds:
\[
\Vert \bigtriangleup\xglm^{[n-1]} \Vert \sim \mathcal{O}(\varepsilon^{p+1}/h), \qquad 
\Vert \bigtriangleup\zglm^{[n-1]} \Vert \sim \mathcal{O}(\varepsilon^{p+1}/h).
\]
From \eqref{eqn:remainder-estimate} we obtain the desired estimates for the remainders $\Vert \Delta_\xglm^{[n]} \Vert$ and $\Vert \Delta_\zglm^{[n]} \Vert$.
\end{proof}

\begin{corollary}[No order reduction]
The global error result \eqref{eqn:spp-order} indicates that the numerical solutions for both $\xx$ and $\zz$ variables converge at the theoretical order $p$:
\[
\xglm^{[n]} - \mathbf{W}\E\, \Nordsieck_p(t_n,h,\xx) \sim \mathcal{O}(h^p), \qquad
\zglm^{[n]} - \mathbf{W}\I\, \Nordsieck_p(t_n,h,\zz) \sim \mathcal{O}(h^p),
\]
for any $\varepsilon \le c h$, $c \le D$. This means that no order reduction is incurred for problems ranging from very stiff ($c \ll 1$) to non-stiff ($c \sim \mathcal{O}(1)$).
\end{corollary}

\begin{remark}[Possible order reduction]
Consider again the iteration \eqref{eqn:bound}. Application of the iteration matrix $K$ times leads to an error amplification factor:
\begin{equation*}
\begin{bmatrix} (\mathbf{V}\E)^K \otimes \Idvar + \mathcal{O}(h) + \mathcal{O}(\varepsilon) & \mathcal{O}\left(\varepsilon\right) \\ 
\mathbf{C}_K
& \prod_{i=1}^K \mathbf{M}\I\left( (h/\varepsilon)\,\gun_\zz\big|_{t_{n+K-i}} \right) + \mathcal{O}(\varepsilon) \end{bmatrix}.
\end{equation*}
The block $\mathbf{C}_K$ is a sum of products of posers of $\mathbf{V}\E$ and $\mathbf{M}\I$.

If we replace assumption \eqref{eqn:eigenvalue-assumption} by the weaker assumption that $\mathbf{M}\I(z)$ is (uniformly) power bounded, and
$\mathbf{C}_K \sim \mathcal{O}(1)$, then the above error amplification matrix is uniformly bounded. Then, by standard arguments, the iterations \eqref{eqn:bound} converge, and give a solution
\[
\Vert \bigtriangleup\xglm^{[n-1]} \Vert \sim \mathcal{O}(\varepsilon^{p+1}/h), \qquad 
\Vert \bigtriangleup\zglm^{[n-1]} \Vert \sim \mathcal{O}(\varepsilon^{p+1}/h^2).
\]
Consider again global error result \eqref{eqn:spp-order}. Under these weaker assumptions, in the nonstiff case $\varepsilon \sim \mathcal{O}(h)$, order reduction in the stiff variable is possible:
\[
\zglm^{[n]} - \mathbf{W}\I\, \Nordsieck_p(t_n,h,\zz) \sim \mathcal{O}(h^{p-1}).
\]
\end{remark}

\subsection{Some useful results}

In this section we collect several useful results used in the proofs.

\begin{lemma}[Eigenvalues and logarihmic norm]
\label{lem:eig-log-norm}
For any matrix $\mathbf{J}$ the real part of its eigenvalues $\lambda_i$ is bounded by the logarithmic norm as follows:
\[
-\mu(\mathbf{-J}) \le   Re\, \lambda_i \le \mu(\mathbf{J}).
\]
\end{lemma}

\begin{theorem}[von Neumann theorem]
\label{thm:vonNeumann}
We restate \cite[Theorem IV.11.2]{Hairer_book_II}. Let $R(z)$ be a rational function bounded on the complex left half plane: 
\[
| R(z) | \le M \quad \forall \,z \in \Co^{-}.
\]
Then, for matrices $\mathbf{J} \in \Re^{\nvar \times \nvar}$:
\[
Re \left( \yy^T\, \mathbf{J}\, \yy \right) \le 0,~ \forall \yy \in \Co^{\nvar}
\qquad \Rightarrow \qquad
\Vert R(\mathbf{J}) \Vert \le M.
\]
\end{theorem}

\begin{theorem}[Matrix-valued von Neumann theorem]
\label{thm:vonNeumann-matrix}
We restate \cite[Theorem V.7.8]{Hairer_book_II}. Let $\mathbf{G} \in \Re^{r \times r}$ positive definite, and let
\[
|\!|\!| Y |\!|\!|_{\mathbf{G}}^2 \coloneqq \sum_{i,j=1}^r \mathbf{G}_{i,j}\, Y_i^T\, Y_j, \quad Y_i \in \Re^{\nvar}, \quad Y = [Y_1^T \dots Y_r^T]^T  \in \Re^{\nvar r}.
\]
Let $\mathbf{M}(z) \in \Re^{r \times r}$ be a matrix whose entries are rational functions of $z$. If
\[
\Vert \mathbf{M}(z) \Vert_{\mathbf{G}} \le M \quad \forall z \in \Co^{-}
\]
then
\[
|\!|\!| \mathbf{M}(\mathbf{J}) |\!|\!|_{\mathbf{G}} \le M \quad \forall \mathbf{J} \in \Re^{\nvar \times \nvar}\,:\, Re \left( \yy^T\, \mathbf{J}\, \yy \right) \le 0,~ \forall \yy \in \Co^{\nvar}.
\]
\end{theorem}

\begin{lemma}[Influence of perturbations on the IMEX-GLM internal stages (Schneider, 1993)]
\label{lem:perturbations-on-internal-stages}
We restate \cite[Theorem 2.2]{Schneider_1993_GLM-SPP}.
Consider the singular perturbation problem \eqref{eqn:SPP}, with an initial solution \eqref{eqn:SPP-epsilon-expansion} that is consistent with respect to \eqref{eqn:dae-eps-0}, \eqref{eqn:expansion-eps-nu}. Apply an IMEX-GLM scheme of order $p$, satisfying the assumptions of Theorems \ref{thm:solution-existence} and \ref{thm:index-k-order} hold, to obtain the stage solutions $X^{[n]},Z^{[n]}$. Apply the same method starting from perturbed initial conditions \eqref{eqn:initial-conditions-assumption}
\begin{equation}
\label{eqn:perturbed-initial-condition}
\widehat{\xglm}^{[n-1]} = \xglm^{[n-1]} + \mathcal{O}(h), \qquad 
\widehat{\zglm}^{[n-1]} = \zglm^{[n-1]} + \mathcal{O}(h),
\end{equation}
and with stage computations affected by perturbations $\delta^{[n]}$ and $\theta^{[n]}$: 
\begin{equation*}
\begin{split}
\widehat{X}^{[n]} &= h\, \left(\mathbf{A}\E \otimes \Idvar \right) \, \fun\big( \widehat{X}^{[n]},\widehat{Z}^{[n]}\big)   + \left(\mathbf{U}\E \otimes \Idvar \right) \, \widehat{\xglm}^{[n-1]} + h\,\delta^{[n]}, \\ 
\varepsilon\,\widehat{Z}^{[n]} &= h\,\left( \mathbf{A}\I \otimes \Idvar \right) \, \gun\big( \widehat{X}^{[n]},\widehat{Z}^{[n]}\big) + \left(\mathbf{U}\I \otimes \Idvar \right) \,\varepsilon\,\widehat{\zglm}^{[n-1]}+ h\,\theta^{[n]}.
\end{split}
\end{equation*}
For $\delta^{[n]} \sim \mathcal{O}(1)$ and $\theta^{[n]} \sim \mathcal{O}(h)$ the difference in the stage solutions is:
\begin{equation*}
\begin{split}
\bigtriangleup X^{[n]} &=   C_1^{[n]}\,\bigtriangleup \xglm^{[n-1]} + \varepsilon\,C_2^{[n]}\, \bigtriangleup \zglm^{[n-1]}  + \mathcal{O}(h   \Vert \delta^{[n]}  \Vert  + h   \Vert \theta^{[n]}  \Vert ), \\ 
\bigtriangleup Z^{[n]} &=   C_3^{[n]}\,\bigtriangleup \xglm^{[n-1]} + \frac{\varepsilon}{h}\, C_4^{[n]}\,\bigtriangleup \zglm^{[n-1]}  + \mathcal{O}(h   \Vert \delta^{[n]}  \Vert  + \Vert \theta^{[n]}  \Vert ),
\end{split}
\end{equation*}
where $C_i^{[n]}$ are families of matrices bounded uniformly for all $n$.
%
%
\end{lemma}

\begin{proof}
The reasoning follows closely the proofs of existence of a numerical solution for implicit Runge-Kutta methods \cite{Hairer_1988_RK-SPP}, \cite[Theorem VI.3.6]{Hairer_book_II} and for implicit GLMs \cite{Schneider_1993_GLM-SPP}. Recall notation \eqref{eqn:initial-conditions} and let
\begin{equation*}
\widehat{\eta}^{[n-1]} \coloneqq \left(\mathbf{U}\E \otimes \Idvar \right) \, \widehat{\xglm}^{[n-1]}, \qquad
\widehat{\zeta}^{[n-1]} \coloneqq \left(\mathbf{U}\I \otimes \Idvar \right) \, \widehat{\zglm}^{[n-1]} . 
\end{equation*}
The nonlinear system \eqref{eqn:existence-nonlinear-system} can be written as a homotopy:
\begin{eqnarray*}
\widehat{X}^{[n]} - \widehat{\eta}^{[n-1]} - h\, \left(\mathbf{A}\E \otimes \Idvar \right) \, \fun\big(\widehat{X}^{[n]},\widehat{Z}^{[n]}\big) 
&=& \tau\,(\widehat{\eta}^{[n-1]} - \eta^{[n-1]} + h\,\delta^{[n]}), \\ 
(\varepsilon/h)\, \bigl( \widehat{Z}^{[n]} - \widehat{\zeta}^{[n-1]} \bigr) - \left( \mathbf{A}\I \otimes \Idvar \right) \, \gun\big(\widehat{X}^{[n]},\widehat{Z}^{[n]}\big) &=& \tau\,\varepsilon\,(\widehat{\zeta}^{[n-1]} - \zeta^{[n-1]} + h\,\theta^{[n]}),
\end{eqnarray*}
where $\tau=0$ gives the base case \eqref{eqn:existence-nonlinear-system} and $\tau=1$ the difference between the perturbed and unperturbed solutions. The proof follows exactly as in \cite{Hairer_1988_RK-SPP}, \cite[Theorem VI.3.6]{Hairer_book_II},\cite{Schneider_1993_GLM-SPP}.
\end{proof}

\begin{lemma}[Influence of perturbations on the IMEX-GLM external stages (Schneider, 1993)]
\label{lem:perturbations-on-external-stages}
We restate part of the proof of \cite[Theorem 2.2]{Schneider_1993_GLM-SPP} as a separate result.
Consider the matrix
\begin{equation*}
\mathbf{M}(t) = \begin{bmatrix} \mathbf{V}\E \otimes \Idvar + \mathcal{O}(h) & \mathcal{O}\left(\varepsilon\right) \\ \mathcal{O}(1)  & \mathbf{M}\I\left( t \right) + \mathcal{O}(\varepsilon) \end{bmatrix}
\end{equation*}
where the spectral radius of the lower right block is uniformly bounded for all times:
\[
\rho\bigl( \mathbf{M}\I\left( t \right) \bigr) \le \alpha < 1 \quad \forall\; t, ~ h/\varepsilon 
\le D.
\]
Then the solution of the following iteration:
\begin{equation*}
\begin{bmatrix} \bigtriangleup\xglm^{[n]} \\ \bigtriangleup\zglm^{[n]} \end{bmatrix}
=  \mathbf{M}(t_n)
\begin{bmatrix} \bigtriangleup\xglm^{[n-1]} \\ \bigtriangleup\zglm^{[n-1]} \end{bmatrix} + 
\begin{bmatrix} \mathcal{O}(\varepsilon^{p+1}) \\ \mathcal{O}(\varepsilon^{p+1}/h) \end{bmatrix},
\qquad \begin{bmatrix} \bigtriangleup\xglm^{[0]} \\ \bigtriangleup\zglm^{[0]} \end{bmatrix}
= \begin{bmatrix} \mathcal{O}(\varepsilon^{p+1}) \\ \mathcal{O}(\varepsilon^{p+1}) \end{bmatrix},
\end{equation*}
has the magnitude
\[
\Vert \bigtriangleup\xglm^{[n-1]} \Vert \sim \mathcal{O}(\varepsilon^{p+1}/h), \quad 
\Vert \bigtriangleup\zglm^{[n-1]} \Vert \sim \mathcal{O}(\varepsilon^{p+1}/h).
\]
\end{lemma}


\section{Discussion}
\label{sec:discussion}

This paper develops new fixed-step convergence results for implicit-explicit general linear methods. We focus on a subclass of schemes that is internally consistent, has high stage order, and favorable stability properties.

The classical convergence analysis reveals that IMEX GLMs in the class of interest converge with the full theoretical order for sufficiently small step sizes. The upper bound for the step size depends on the method coefficients and on non-stiff Lipschitz constants, but is independent of the stiffness of component $\fun\I$. 

Convergence analysis for IMEX-GLMs applied to index-1 differential algebraic problems reveals that the methods in the class of interest converge with the full theoretical order. An order reduction (to stage order) for the algebraic variable is possible when the implicit stability matrix at infinity has (simple) eigenvalues of magnitude one.

IMEX GLMs applied to singular perturbation problems have a unique solution under general assumptions, and methods in the class of interest converge with the full theoretical order for both the stiff and the non-stiff variables.

\section*{Acknowledgements}

This work has been supported in part by NSF through awards NSF ACI–1709727 and NSF CCF–1613905, and by the Computational Science Laboratory at Virginia Tech.

\bibliographystyle{elsarticle-num}

\end{document}